\documentclass[final,3p,times]{elsarticle}

\usepackage{amssymb}
\usepackage{amsmath}
\usepackage{amsthm}
\usepackage{mathtools}

\usepackage{hyperref}
\usepackage{lipsum}
\usepackage{amsfonts}
\usepackage{graphicx}
\usepackage{epstopdf}
\usepackage{algorithmic}
\usepackage{multirow}
\usepackage{multicol}
\usepackage{caption}
\usepackage{subcaption}
\usepackage{setspace}
\usepackage{xcolor}
\usepackage{cleveref}
\usepackage{bm}
\usepackage{algorithm}
\usepackage{algorithmic}
\usepackage{pgfplots}
\pgfplotsset{compat=1.18}
\usepackage{overpic}

\usepackage[normalem]{ulem}

\newtheorem{theorem}{Theorem}[section]   

\newtheorem{remark}[theorem]{Remark}

\ifpdf
  \DeclareGraphicsExtensions{.eps,.pdf,.png,.jpg}
\else
  \DeclareGraphicsExtensions{.eps}
\fi

\renewcommand{\epsilon}{\varepsilon} 
\renewcommand{\S}{\mathcal{S}}

\renewcommand{\phi}{\varphi}

\renewcommand{\O}{\mathcal{O}}

\newcommand{\Z}{\mathbb{Z}}
\newcommand{\R}{\mathbb{R}}

\newcommand{\N}{\mathbb{N}}

\newmuskip\pFqmuskip
\newcommand*\pFq[6][8]{%
  \begingroup 
  \pFqmuskip=#1mu\relax
  \mathcode`\,=\string"8000
  \begingroup\lccode`\~=`\,
  \lowercase{\endgroup\let~}\pFqcomma
  {}_{#2}F_{#3}{\left[\genfrac..{0pt}{}{#4}{#5};#6\right]}%
  \endgroup
}
\newcommand{\pFqcomma}{\mskip\pFqmuskip}

\newcommand{\revision}[1]{\textcolor{black}{#1}}
\newcommand{\delete}[1]{}

\allowdisplaybreaks

\journal{Journal of Computational Physics}

\begin{document}

\begin{frontmatter}



\title{Fully Explicit Stabilized Runge--Kutta Methods for Allen--Cahn and Cahn--Hilliard Dynamics}



\author[label1]{Zheng Tan}
\ead{zhengtan@math.ucla.edu}

\author[label2]{Tariq D. Aslam}
\ead{aslam@lanl.gov}

\author[label1,label3]{Andrea L. Bertozzi\corref{cor1}}
\ead{bertozzi@ucla.edu}

\cortext[cor1]{Corresponding author}

\affiliation[label1]{organization={Department of Mathematics, University of California, Los Angeles},
            state={CA},
            postcode={90095}, 
            country={USA}}
\affiliation[label2]{organization={Theoretical Division, Los Alamos National Laboratory, Los Alamos},
            state={NM},
            postcode={87545}, 
            country={USA}}
\affiliation[label3]{organization={Department of Mechanical and Aerospace Engineering, University of California, Los Angeles},
            state={CA},
            postcode={90095}, 
            country={USA}}

\begin{abstract}
Stiff phase-field dynamics are typically solved by semi-implicit or fully implicit time
discretizations, since the interfacial scales in the Allen--Cahn and Cahn--Hilliard equations
impose severe step-size restrictions on standard explicit methods. In this paper, we investigate a
different strategy based on stabilized Runge--Kutta solvers, in particular the Runge–Kutta–Legendre (RKL) and Runge–Kutta–Gegenbauer (RKG) families. These schemes possess enlarged stability regions and therefore provide a natural explicit framework for diffusion-dominated gradient flows without requiring implicit solves. For the Allen--Cahn equation, combined with Strang splitting, we establish a maximum principle in 1-dimensional case, and prove energy dissipation in all the dimensions. For the modified energy dissipation analysis, we further introduce two new first-order stabilized schemes of RKL and RKG type and derive a novel eigenvalue estimate for the second-order RKG method. We also conduct extensive numerical experiments with both constant and adaptive time stepping, on regular domains as well as on an L-shaped irregular domain. The numerical results show that stabilized explicit Runge--Kutta methods
provide an effective and flexible approach for the numerical simulation of Allen--Cahn and Cahn--Hilliard dynamics while preserving key structural properties. The Python implementation is available at \url{https://github.com/ZT220501/SRK-for-Gradient-Flows}.
\end{abstract}

\begin{graphicalabstract}
\end{graphicalabstract}

\begin{highlights}
\item Stabilized Runge--Kutta methods for fully explicit, non-spectral solution approach of the Allen--Cahn and the Cahn--Hilliard Dynamics.
\item The Maximum principle is proved when applying to the Allen--Cahn equation.
\item Modified energy dissipation is proved for two Strang splitting orders.
\item The method applies to irregular domains without Fourier spectral methods.
\end{highlights}

\begin{keyword}
Stabilized Runge--Kutta \sep Gradient flows \sep Energy stability \sep Maximum principle preserving method

\end{keyword}

\end{frontmatter}



\section{Introduction}
\label{intro}
Phase-field models are prevalent in science and engineering, particularly in material science and fluid dynamics, see for example~\cite{AllenCahn,CahnHilliard,annurev:/content/journals/10.1146/annurev.matsci.32.112001.132041,Golubovic_Levandovsky_Moldovan_2011,doi:10.1142/S0218202596000341, PENROSE199044, GOMEZ2013217}. Therefore, it is important to efficiently and accurately solve the equations numerically.
In this work, we consider two of the most well-known phase-field equations, namely the Allen--Cahn equation
\begin{equation}
\label{eqn:AllenCahn}
\frac{\partial u}{\partial t}=\epsilon^2\nabla^2u-F'(u), \quad x\in\Omega, t\in[0, T],
\end{equation}
and the Cahn--Hilliard equation
\begin{equation}
\label{eqn:CahnHilliard}
\frac{\partial u}{\partial t}=\nabla^2\left(-\epsilon^2\nabla^2u+F'(u)\right),\quad x\in\Omega, t\in[0, T],
\end{equation}
with initial condition $u(x, 0)=u_0(x)$. Here $\Omega\subset\R^n$ is a bounded domain and $T>0$ is a finite time. Throughout this paper, we use the classical double-well potential $F(u):=\frac{1}{4}(u^2-1)^2$. We use either the periodic boundary conditions, or the homogeneous Neumann boundary condition
$$\frac{\partial u}{\partial\bm{n}}\Big|_\Omega=\frac{\partial\mu}{\partial\bm{n}}\Big|_\Omega=0,$$
where $\mu=-\epsilon^2\nabla^2u+F'(u)$ and $\bm{n}$ is the outward normal of $\partial\Omega$.

The Allen--Cahn and the Cahn--Hilliard equations are well-known to be the gradient flows of the Ginzburg--Landau energy 
\begin{equation}
\label{eqn:ginzburg_landau_energy}
E(u(t))=\int_\Omega\frac{\epsilon^2}{2}|\nabla u|^2+F(u)dx.
\end{equation}
More precisely, the $L^2$-gradient flow results in the Allen--Cahn equation, while the $H^{-1}$-gradient flow results in the Cahn--Hilliard equation. Therefore, a successful numerical scheme should also preserve the energy dissipation property to obtain physically correct solutions; more precisely, the numerical solution should satisfy $$E(u(t_1))\le E(u(t_0)),\quad\forall t_1\ge t_0.$$ 
We call a scheme \textit{energy stable} if the above condition is satisfied. Another important property of the second-order parabolic Allen--Cahn equation is the maximum principle: if $\|u(\bm{x}, 0)\|_\infty\le1$, then $\|u(\bm{x}, t)\|_\infty\le1$ for all $t>0$. Therefore, it is also essential for a numerical scheme to satisfy the maximum principle for the Allen--Cahn equation.

\delete{Since gradient flows involve harmonic or even biharmonic operators, along with strong nonlinearity, it is difficult to design a numerical scheme for their simulations while preserving the energy stability and/or maximum principle.} \revision{Preserving numerical energy stability is central to schemes for gradient flow simulations.} An abundance of numerical evidence shows that spurious oscillations may occur when the energy stability condition is violated. \revision{Numerical effort has been devoted for energy stable schemes development, including~\cite{ConvexSplitting, IMEX, XU201936, SAV}}. Since the spatial derivative is second or fourth order, the stability constraint limits the usage of classical explicit methods, leading to the need to solve an implicit system at each step. The recently deployed exponential time differencing Runge--Kutta schemes (ETDRK)~\cite{ETDRKEnergyStable, ETDRK2MaximumPrinciple, ETDRKEnergyStableHigherOrder}, although they avoid implicit solves and monotonically dissipate energy, rely heavily on the Fourier spectral methods so that on non-periodic irregular domains, it degenerates into an implicit method because of the need for a stabilization parameter.

In this paper, we explore the potential of fully explicit stabilized Runge--Kutta methods called Runge--Kutta--Legendre (RKL) and Runge--Kutta--Gegenbauer (RKG) methods for solving the Allen--Cahn and Cahn--Hilliard equations. The RKL and RKG methods are super-time-stepping multistage methods that admit stability constraints to be relaxed quadratically with their number of stages. We prove that the RKL and RKG methods preserve the maximum principle for the Allen--Cahn equation when combined with the Strang splitting, under periodic or Neumann boundary conditions. Modified energy dissipations for the RKL and RKG methods are proved when they are applied to the Allen--Cahn equation, and the dissipated energy is shown to converge to the exact energy. We also conduct numerous experiments to demonstrate the behaviors when the RKL and RKG methods are applied to the gradient flows. Specifically, an irregular domain experiment is conducted to demonstrate the effectiveness of the methods for a case where the Fourier spectral method cannot be easily applied.

This paper is organized as follows. In~\Cref{sec:literature_review} we review some of the well-known numerical schemes that efficiently and accurately solve the gradient flows while preserving the energy stability, with their advantages and limitations.~\Cref{sec:SRK} introduces the RKL and RKG methods and describes how we apply them to the gradient flows. Following that, in~\Cref{sec:theorems} we prove the maximum principle and the modified energy dissipation for RKL and RKG schemes. We conduct comprehensive experiments in~\Cref{sec:experiment} to demonstrate the behavior of RKL and RKG methods when applied to gradient flows. Concluding remarks and potential research directions are presented in~\Cref{sec:conclusion}.

\section{Review of schemes for simulating gradient flows}
\label{sec:literature_review}
Various schemes have been designed in order to satisfy the maximum principle of the Allen--Cahn equation and energy stability of the gradient flows. In this section, we present several of the most representative methods, in roughly chronological order.

\subsection{Convex splitting}
\label{subsec:convex_splitting}
One of the pioneering works for energy dissipating Cahn--Hilliard simulation is the convex splitting technique~\cite{ConvexSplitting}. Assume that the energy $F(u)$ can be split into the difference of two convex functions; more precisely, $F(u)=F_1(u)-F_2(u)$ where both $F_1, F_2$ are convex. Then one treats the positive convex part to be implicit and negative convex part to be explicit, so that the first-order convex splitting method for Cahn--Hilliard reads
\begin{equation}
\label{eqn:convexsplitting}
\frac{u^{n+1}-u^n}{\tau}=\nabla^2\mu^{n+1};\quad\mu^{n+1}=\left(-\epsilon^2\nabla^2u^{n+1}+F_1'(u^{n+1})-F_2'(u^n)\right).
\end{equation}
It \delete{was}\revision{is} shown that~\cite{ConvexSplitting}
$$E(u^{n+1})-E(u^n)\le-\tau\|\nabla\mu^{n+1}\|^2,$$
so that the scheme dissipates the energy as desired. However, a nonlinear system needs to be solved at every step, leading to inefficient computation.

\subsection{Implicit-Explicit (IMEX) methods}
Another approach, to some extent similar to the convex splitting idea, is the implicit-explicit (IMEX) method~\cite{IMEX}. Instead of splitting the energy $F(u)$, the IMEX approach directly treats the linear part implicitly while treating the nonlinear part explicitly. Therefore, only a linear system needs to be solved at each step.

Several works have studied IMEX methods for parabolic equations and gradient flows~\cite{IMEX,vollmayrlee2003fast,IMEXSixth,FuTangYang2024}. Although only linear system solving is required, the IMEX method may still only attain conditional energy dissipation, unless additional stabilization techniques~\cite{JieShenXiaofengYang2010} are used. Therefore, more unconditional energy dissipation techniques are developed, as will be introduced in the following subsections.


\subsection{Lagrange multiplier approaches}
One of the novel approaches for solving gradient flows utilizes the idea of Lagrange multiplier. X. Yang and his collaborators \delete{applies}\revision{apply} the Lagrange multiplier approach to the class of gradient flows with free energies bounded below~\cite{XU201936,YANG2016294,doi:10.1142/S0218202521500184,YANG2017691,YANG20171005,YANG201880,YANG2017104,yang2017numerical,YU2017665} and \delete{developed}\revision{develop} the invariant energy quadratization (IEQ) method. To illustrate the approach, consider the Cahn-Hilliard equation written in the form of 
\begin{equation}
\label{eqn:cahn-hilliard-ieq}
\frac{\partial u}{\partial t}=\nabla^2\mu;\quad\mu=-\nabla^2u+F'(u).
\end{equation}
With the free energy $F(u)$ uniformly bounded below such that $F(u)\ge-C_0\in\R$, one introduces the Lagrange multiplier $q(t,x;u)=\sqrt{F(u)+C_0}$, and rewrites~\Cref{eqn:cahn-hilliard-ieq} as
\begin{equation}
\label{eqn:IEQ}
\frac{\partial u}{\partial t}=\nabla^2\mu;\quad\mu=-\nabla^2u+\frac{q}{\sqrt{F(u)+C_0}}F'(u);\quad\frac{\partial q}{\partial t}=\frac{F'(u)}{2\sqrt{F(u)+C_0}}\frac{\partial u}{\partial t}.
\end{equation}
The IEQ method solves~\eqref{eqn:IEQ} with numerical methods such as Crank-Nicolson or BDF schemes. Taking the inner product of the above with $\mu, u_t$, and $2q$ respectively, the IEQ method naturally satisfies the modified energy dissipation law
$$\frac{d}{dt}\Big(\frac{1}{2}\|\nabla u\|^2+\int_\Omega q^2dV\Big)=-\|\nabla\mu\|^2,$$
so that the energy stability is automatically guaranteed.

Though unconditionally energy stable, the IEQ method requires solving a variable coefficient system of equations every step. This disadvantage is overcome by the scalar auxiliary variable (SAV) approach~\cite{SAV}. Instead of assuming the free energy $F(u)$ is bounded below, the SAV approach assumes that the nonlinear part of energy $E_1(u):=\int_\Omega F(u)dV$ is bounded from below such that $E_1(u)\ge-C_0\in\R$. By introducing the Lagrange multiplier $r(t;u)=\sqrt{E_1(u)+C_0}$, the Cahn--Hilliard equation~\eqref{eqn:cahn-hilliard-ieq} can be rewritten as 
\begin{equation}
\label{eqn:SAV}
\frac{\partial u}{\partial t}=\nabla^2\mu;\quad\mu=-\nabla^2u+\frac{r}{\sqrt{E_1(u)+C_0}}F'(u);\quad\frac{\partial r}{\partial t}=\frac{1}{2\sqrt{E_1(u)+C_0}}\int_\Omega F'(u)\frac{\partial u}{\partial t}dV.
\end{equation}
Using Crank-Nicolson or BDF methods, only scalar-coefficient systems are needed to be solved at each step for the SAV method. A similar modified energy dissipation law also holds for the SAV approach
$$\frac{d}{dt}\Big(\frac{1}{2}\|\nabla u\|^2+r^2(t)\Big)=-\|\nabla\mu\|^2.$$
Extensive follow-up works of the SAV method~\cite{Bai2023GinzburgHminus1SAV,CHENG2020113070,GSAV,doi:10.1137/24M1659303,Di2023VariableIMEXBDF2SAV,HUANG2023112034,LI202362,MSAV,Liu2024EOPSAV,SAVReview,RGSAV} \delete{demonstrated}\revision{demonstrate} its effectiveness across various gradient flow problems. However, the method still requires solving an implicit system of equations at each step since explicit solvers cannot be applied to~\eqref{eqn:IEQ} or~\eqref{eqn:SAV} directly, leading to the need for solving two Poisson-type equations per step~\cite{SAV}.

\subsection{Exponential time differencing (ETD) methods}
A popular and active line of research is the use of exponential time differencing (ETD) methods~\cite{ExponentialIntegrators} for solving the phase-field equation. The phase-field equation has the semi-linear right hand side structure
$$\frac{\partial u}{\partial t}=G(Lu+f(u)),$$
where $G$ is a non-positive operator, $L$ is a positive-definite operator, and $f$ is a (nonlinear) function. Using Duhamel's principle, the \textit{exact} solution is
\begin{equation}
\label{eqn:duhamel}
u(x, t)=e^{GL(t-t_0)}u(x_0, t_0)-e^{GL(t-t_0)}\int_{t_0}^te^{-GL(s-t_0)}Gf(u(x, s))ds.
\end{equation}
The ETD methods are then developed by approximating the integral. A major advantage of the ETD methods is that the linear stiff part is exactly integrated, so that larger timesteps can be chosen after the matrix exponentials are computed using the Fourier transform.

The most widely used ETD methods are the ETD Runge--Kutta (ETDRK) methods, which use RK schemes to approximate the integral in~\eqref{eqn:duhamel}. Much of the existing literature has studied the ETDRK methods for semi-linear parabolic equations and gradient flows~\cite{ETDRK2MaximumPrinciple,CoxMatthews,ETDRKSemilinear,ETDRKEnergyStable,ETDRKEnergyStableHigherOrder,ETDRKMatrixValued}. The discrete maximum principle for first- and second-order ETDRK methods (ETDRK1, ETDRK2) \delete{were}\revision{are} proved for the second-order parabolic equations~\cite{ETDRKSemilinear}. The dissipation of the original energy for the Ginzburg--Landau gradient flow was proved for ETDRK1, ETDRK2~\cite{ETDRKEnergyStable}; later a general sufficient criterion for ETDRK of arbitrary order \delete{was proposed}\revision{is later proposed}, along with a novel energy dissipating ETDRK3 method~\cite{ETDRKEnergyStableHigherOrder}. The existence of energy dissipating ETDRK4 or higher order ETDRK methods is still an open problem.

Although implicit equation solvers are not needed for the ETDRK methods, the success of ETDRK methods for gradient flows largely relies on the effectiveness in computing the matrix exponentials under the periodic boundary condition within a regular domain. Since the Laplacian operator is diagonalizable in the Fourier domain, the computational cost is only $\O(N^2\log N)$ if there are $N$ grids for each direction in two dimension~\cite{ETDRK2MaximumPrinciple}. However, if a non-periodic boundary condition is applied, the Fourier spectral method cannot be directly applied~\cite{spectralmethodtextbook}. Moreover, a stabilization parameter is needed in practice in order for the method to behave well and preserve the maximum principle~\cite{ETDRKEnergyStable}, making the method essentially implicit with the need for tuning the stabilization parameter.

\section{Stabilized Runge--Kutta methods}
\label{sec:SRK}
Due to the stiff nature of gradient flows, the existing methods in~\Cref{sec:literature_review} either use (semi-)implicit approaches, or utilize the exponential time differencing approach along with the Fourier transform and stabilization technique. Although linear and energy stability can be guaranteed with mild or no constraints, the (semi-)implicit methods require solving a linear or nonlinear system at each step, while the exponential time differencing approach requires computing matrix exponential, leading to heavy computational overhead.

In this section, we utilize a class of stabilized Runge--Kutta (SRK) methods for solving the Allen--Cahn and Cahn--Hilliard equations, more specifically the Runge--Kutta--Legendre (RKL) and Runge--Kutta--Gegenbauer (RKG) methods. The SRK methods are \textit{fully explicit} methods, and do not rely on any computations in the Fourier domain or stabilization techniques. Although the SRK methods are only conditionally stable due to their explicit nature, the region of absolute stability of the SRK methods grow proportionally to $s^2$ where $s$ is the number of sub-stages, leading to mild stability constraints in practice. Rigorous analysis shows that along with Strang splitting, the RKL and RKG methods preserve the maximum principle under the linear stability constraints. Following similar arguments as in~\cite{StrangSplittingAC, GAO2026116868}, modified energy dissipation for two different ways of Strang splittings are proved for the Allen--Cahn equation, with new first-order methods and an estimation for the second-order RKG (RKG2) method. The modified energy is shown to converge to the exact energy as $\tau\to0$.

We note that the existing literature has briefly explored the use of super-time-stepping method for Cahn--Hilliard equation simulation~\cite{CHENG201545}. However, they \delete{utilized}\revision{utilize} Runge--Kutta--Chebyshev method and the DUMKA3 methods and only simulated the fourth-order Cahn--Hilliard and molecular beam epitaxy (MBE) equations; such methods have been numerically shown to violate the discrete maximum principle~\cite{RKL2}, so that they cannot be extended to the Allen--Cahn equation simulation. Moreover, our method is independent of spatial discretization, so that it can be applied to more generally shaped domains in addition to rectangular periodic domains, as will be shown in~\Cref{sec:experiment}. Finally, we provide an energy dissipation theoretical guarantee for the Allen--Cahn equation simulation with the SRK, for different ways of Strang splittings.

\subsection{Stabilized Runge--Kutta methods overview}
\label{subsec:rkl_rkg}

The Runge--Kutta--Legendre methods~\cite{RKL2} and Runge--Kutta--Gegenbauer methods~\cite{OSULLIVAN_RKG2019209, SKARAS_RKG2021109879} are super-time-stepping methods originally designed to solve for an ordinary differential equation (ODE) system 
\begin{equation}
\label{autoode}
    \frac{du}{dt}=\bm{M}u,
\end{equation}
where $\bm{M}$ is a symmetric, constant coefficient matrix that represents the discretization of a parabolic operator; the operator $\bm{M}$ is required to have non-positive, real eigenvalues to ensure convergence. It \delete{was}\revision{is} also shown~\cite{RKL2, OSULLIVAN_RKG2019209, SKARAS_RKG2021109879} that the RKL and RKG methods are effective for variable coefficient problems, and later~\cite{tan2025singularity, tan2026stork} apply this class of method for broader applications.

We present the derivation of the RKL and RKG methods. The RKL and RKG methods are based on the shifted Legendre and Gegenbauer polynomials, respectively. The general form of the RKL and RKG methods is
\begin{align}
u^{t+\tau}=R_s(\tau\bm{M})u^t,
\end{align}
where $R_s(\tau\bm{M})$ is the ``stability polynomial" of the numerical method. Numerical stability is ensured if $|R_s(\tau\lambda)|\le1$ for all eigenvalues $\lambda$ of $M$; therefore, boundedness of Legendre polynomials and Gegenbauer polynomials enables them to be the stability polynomial choices~\cite{RKL2, SKARAS_RKG2021109879}. By matching the leading order terms in Taylor expansion as in~\cite{RKL2, SKARAS_RKG2021109879}, the RKL1, RKL2, RKG1, and RKG2 methods, respectively, take the form of
\begin{align}
\label{eqn:RKL1}
u^{t+\tau}=P_{s}\Big(I+\frac{2\tau}{s^2+s}\bm{M}\Big)u^t;
\end{align}
\begin{align}
\label{eqn:RKL2}
u^{t+\tau}=\Big[\Big(1-\frac{s^2+s-2}{2s(s+1)}\Big)I+\frac{s^2+s-2}{2s(s+1)}P_{s}\Big(I+\frac{4\tau}{s^2+s-2}\bm{M}\Big)\Big]u^t;
\end{align}
\begin{align}
\label{eqn:RKG1}
u^{t+\tau}=\frac{2}{(s+1)(s+2)}C^{(\lambda)}_{s}\Big(I+\frac{4\tau}{s^2+3s}\bm{M}\Big)u^t;
\end{align}
\begin{align}
\label{eqn:RKG2}
u^{t+\tau}=\Big[\Big(1-\frac{2(s-1)(s+4)}{3s(s+3)}\Big)I+\frac{4(s-1)(s+4)}{3s(s+1)(s+2)(s+3)}C^{(\lambda)}_{s}\Big(I+\frac{6\tau}{(s+4)(s-1)}\bm{M}\Big)\Big]u^t,
\end{align} 
where $P_s(x)$ is the $s$-th degree Legendre polynomial and $C_s^{(\lambda)}(x)$ is the $s$-th degree Gegenbauer polynomial with parameter $\lambda$ defined in the standard literature. In all of our experiments, $\lambda$ is chosen to be $\frac{3}{2}$ in order to be consistent with~\cite{SKARAS_RKG2021109879}.

To obtain a Runge-Kutta formulation of the RKL and RKG methods,~\cite{RKL2, SKARAS_RKG2021109879} \delete{used}\revision{use} the recurrence relationship of the Legendre and Gegenbauer polynomials. For example, the RKL2 method \delete{used}\revision{uses} the Legendre polynomial recurrence relationship
\begin{equation}
    sP_s(x)=(2s-1)xP_{s-1}(x)-(s-1)P_{s-2}(x).
\end{equation}
Together with the initial Legendre polynomials, after rearranging terms, one superstep of the $s$-stage RKL2 method can be written as
\begin{align}
\label{eqn:RK-presentation}
\begin{split}
    &Y_0=u(t_0), \quad Y_1=Y_0+\Tilde{\mu}_1\tau\bm{M}Y_0,\\
    &Y_j=\mu_jY_{j-1}+\nu_jY_{j-2}+(1-\mu_j-\nu_j)Y_0+\Tilde{\mu}_j\tau\bm{M}Y_{j-1}+\Tilde{\gamma}_j\tau\bm{M}Y_0, 2\le j\le s,\\
    &u(t_0+\tau)=Y_s,
\end{split}
\end{align}
where the parameters are specified in~\cite{RKL2}. Similar derivations work for the RKL1 method~\cite{RKL2} and the RKG methods~\cite{SKARAS_RKG2021109879} using the recurrence relation between the Legendre and Gegenbauer polynomials of different degrees.

Although the formulas appear complicated, the method is entirely explicit and can be easily implemented with only three stages of the solution storage required. Moreover, the stability constraint is much more relaxed than a usual explicit method: denote $\tau_{expl}:=\frac{2}{\lambda_{\max}}$ where $\lambda_{\max}$ is the maximum eigenvalue of the operator $\bm{M}$; as shown in~\cite{RKL2}, the largest possible timestep is $\tau_{\max}=\tau_{expl}\frac{s^2+s-2}{4}$. The restriction is relaxed quadratically with respect to $s$, so that when we choose $s$ large enough, the restriction on $\tau$ related to $dx$ is almost negligible. Moreover, it is simple to control the stage number $s$ in the implementation. These advantages lead to the adaptive use of $s$ in our complete algorithm.

Despite the fact that RKL and RKG methods were primarily designed to ensure the discrete maximum principle under different boundary conditions, the theoretical proof of their discrete maximum principle remained unsolved for a long time. Our recent work~\cite{tan2025singularity} fully proved the discrete maximum principle of RKL1, RKL2, RKG1, and RKG2 methods when applied to the heat equation, under the periodic or homogeneous Neumann boundary condition for rectangular domains. The discrete maximum principle of the Allen--Cahn equation therefore follows, as shown in~\Cref{subsec: max_principle}.

\subsection{Stabilized Runge--Kutta methods applied to gradient flows}
With the above description, we now apply the stabilized Runge--Kutta methods to the Allen--Cahn and Cahn--Hilliard equations. For the Allen--Cahn equation~\eqref{eqn:AllenCahn}, since the nonlinear part is a double-well potential, we utilize the Strang splitting~\cite{strangsplitting} technique so that the nonlinear part is solved exactly, while the numerical monotonicity can be guaranteed easily as will be shown in~\Cref{thm:AllenCahnMaximumPrinciple}. By writing 
$$\S_L(\tau)(u^n)=SRK(u^n),\quad \S_N(\tau)(u^n):=\frac{u^n}{\sqrt{(1-(u^n)^2)e^{-2\tau}+(u^n)^2}},$$
our Strang splitting scheme reads
$$u^{n+1}=\S_N(\tau/2)\S_L(\tau)\S_N(\tau/2)u^n,$$
with the discrete Laplacian used being the standard central difference operator.

For the Cahn--Hilliard equation, since \revision{in this paper} we in general do not want to use the Fourier spectral method for spatial discretization, we choose to not use Strang splitting as in~\cite{CHENG201545} but instead discretize the spatial domain with standard central difference, and apply the SRK method with $Mu$ corresponding to the discretized right-hand side of the Cahn--Hilliard equation. 

\section{Maximum Principle and Modified Energy Dissipation}
\label{sec:theorems}
In this section, we prove the maximum principle and the energy dissipation when applying RKL and RKG methods to the Allen--Cahn equation via Strang splitting.

\subsection{Maximum principle for Allen--Cahn via Strang splitting}
\label{subsec: max_principle}
In this section, we present the proof of the discrete maximum principle for the RKL and RKG methods when applied to the Allen--Cahn equation. More precisely, we want to show that if the initial condition is uniformly bounded by 1, then $\|u(t)\|_\infty\le1$ for any further time $t>0$. The discrete analogue is therefore described and proved in the following theorem, utilizing the power of Strang splitting~\cite{strangsplitting} and the monotone stability results proved in~\cite{tan2025singularity}.

\begin{theorem}
\label{thm:AllenCahnMaximumPrinciple}
Consider the \revision{1D} Allen-Cahn equation with periodic or homogeneous Neumann boundary conditions
$$\frac{\partial u}{\partial t}=\epsilon^2\bm{L}u-F'(u),$$
where $F(u)=\frac{1}{4}(u^2-1)^2$, and $\bm{L}$ is the central difference approximation to the second-order spatial derivative. 
\\
Consider the Strang splitting along with any SRK method when solving the above equation
$u^{n+1}=\S_L(\tau/2)\S_N(\tau)\S_L(\tau/2)u^n$ or $u^{n+1}=\S_N(\tau/2)\S_L(\tau)\S_N(\tau/2)u^n$. 
Assume $\|u^0\|_\infty\le1$. Then $\|u^n\|_\infty\le 1$ under the linear stability constraint of the SRK method for all $n\in\N$.
\end{theorem}
\begin{proof}
Notice that for each spatial grid $i$, $1-(u_i^k)^2\ge0$ by assumption, therefore
$$|S_N(\tau)u_{i}^n|=\left|\frac{u_i^n}{\sqrt{(1-(u_i^n)^2)e^{-2\tau}+(u_i^n)^2}}\right|\le\frac{|u_i^n|}{\sqrt{(u_i^n)^2}}=1,$$
and therefore the nonlinear update $S_N(\tau)$ preserves the maximum principle for any $\tau>0$. \delete{Now using the fact that second-order central difference approximation in higher dimensions is the sum of second-order central difference approximations in each dimension, by}\revision{By} Theorem 5.7 and Corollary 5.12 in~\cite{tan2025singularity}, we know that the linear update $\S_L(\tau)$ also preserves the maximum principle, therefore the maximum principle of the Strang splitting follows.
\end{proof}

\begin{remark}
The theorem only works for the most commonly used periodic and homogeneous Neumann boundary conditions. It is known that the RKL methods are not monotonically stable for the homogeneous Dirichlet boundary condition, and therefore cannot be guaranteed to preserve the maximum principle of the Allen--Cahn equation. For the RKG methods, it \delete{was}\revision{is} conjectured in~\cite{SKARAS_RKG2021109879} that they are monotonically stable under the homogeneous Dirichlet boundary condition, while no theoretical proof has been established yet.
\end{remark}

\begin{remark}
The above theorem only works for 1D because of the limitation of the theorem in~\cite{tan2025singularity}. Nevertheless, numerical evidence suggests that the maximum principle also holds a for 2D regular \revision{domain with periodic boundary condition} and irregular domains under the Neumann boundary condition. \delete{and we conjecture that in fact both RKL and RKG should satisfy the maximum principle under the Neumann boundary condition in general domains and RKG should also satisfy the maximum principle in addition under the Dirichlet boundary condition.} \revision{However, direct computations reveal that RKL2 does not always preserve numerical monotonicity in 2D even with the CFL condition.s
}
\end{remark}

\subsection{Modified energy dissipation for the Allen--Cahn equation}
\label{subsec:energy_dissipation}
We now prove the modified energy dissipation of the Allen--Cahn equation. For the linear-nonlinear-linear Strang splitting, we are inspired by the work ~\cite{StrangSplittingAC, GAO2026116868}\delete{; however we need to modify the order of the Strang splitting with an inexact solution of the linear part}\revision{, with an inexact solution of the linear part}. Since the computational cost of the linear part often dominates the cost of the nonlinear part, we propose a different modified energy for the more practical splitting $u^{n+1}=\S_N(\tau/2)\S_L(\tau)\S_N(\tau/2)u^n$ and prove that it dissipates and converges to the exact energy. New first-order RKL and RKG methods are proposed in order to fit into the proof framework, along with an estimation of the RKG2 method spectrum. Moreover, throughout this section we only require the discrete finite difference Laplacian operator $L$ to be real, symmetric, and negative-definite; the domain is not required to be regular.

\subsubsection{Linear-Nonlinear-Linear Strang splitting}
We begin by considering the linear-nonlinear-linear splitting $u^{n+1}=\S_L(\tau/2)\S_N(\tau)\S_L(\tau/2)u^n$. 
\begin{theorem}[Modified Energy Dissipation, linear-nonlinear-linear]
Define the modified energy as
\begin{equation}
\label{eqn:lnl_energy}
\Tilde{E}(u^n):=\frac{1}{2\tau}\Big\langle(I-S^2_L(\tau/2))u^n,u^n\Big\rangle+\langle g(S_L(\tau/2)u^n),1\rangle,
\end{equation}
where
\begin{equation}
\label{eqn:modified_nonlinear_energy}
g(x):=\frac{1}{4}+\frac{1}{2\tau}x^2-\frac{e^\tau}{\tau(e^{2\tau}-1)}\Big(\sqrt{1+(e^{2\tau}-1)x^2}-1\Big).
\end{equation}
Then under the linear-nonlinear-linear Strang splitting defined above, we have
\begin{enumerate}
    \item $\Tilde{E}(u^{n+1})\le\Tilde{E}(u^{n})$;
    \item For fixed grid size $h$, $\Tilde{E}(u^n)$ converges to $E(u^n)$ uniformly as $\tau\to0$; more precisely, for $0\le \tau\le\min\{1, \text{CFL constraint}\}$, assuming enough regularity of $u^0$, we have that
$$\underset{n\ge0}{\sup}|\Tilde{E}(u^n)-E(u^n)|\lesssim \tau.$$
\end{enumerate}
\end{theorem}
\begin{proof}
We will first prove that $\Tilde{E}(\Tilde{u}^{n+1})\le\Tilde{E}(\Tilde{u}^{n})$. By construction, the spectrum of $S_L(\tau/2)$ is within $[-1, 1]$; therefore, we have that $I-S_L^2(\tau/2)$ is positive semi-definite. The non-negativity of the energy and its dissipation result then directly follows from the proof of Theorem 2.1 in~\cite{GAO2026116868}.
To show the uniform convergence, we follow the procedure as in~\cite{StrangSplittingAC}. Define $\Tilde{u}^n:=\S_L(\tau/2)u^n$; taking the difference between $\Tilde{u}^n$ and $u^n$ and using Taylor expansion,
\begin{align*}
\|\Tilde{u}^n-u^n\|=&\|\S_L(\tau/2)u^n-u^n\|=\left\|\left(I+\frac{\tau}{2}L+\mathcal{O}(\tau^2)\right)u^n-u^n\right\|\lesssim\,\mathcal{O}(\tau).
\end{align*}
Notice that $\|\Tilde{u}^n\|_\infty\le\|S_L(\tau/2)\|_\infty\|u^n\|_\infty\le C\|S_L(\tau/2)\|_2\le C$ for some grid size dependent constant $C$ uniformly when the CFL condition is satisfied. \delete{By the same derivation as in~\cite{StrangSplittingAC}, and using the fact that}\revision{By direct calculation, we have that}
$$g(\Tilde{u}^n)=\frac{1}{4}-\frac{1}{2}(\Tilde{u}^n)^2+\frac{1}{4}(\Tilde{u}^n)^4+\O(\tau),$$
we can bound the nonlinear part of the energy \delete{using the same strategy as in~\cite{StrangSplittingAC}.}\revision{using the fact that
\begin{align*}
\left|\int_{\mathbb{T}^d}\Tilde{F}(\Tilde{u})dx-\int_{\mathbb{T}^d}\frac{1}{4}((u^n)^2-1)^2dx\right|\lesssim\tau.
\end{align*}
}
\\
For the linear part, we can bound the difference by
\begin{align*}
&\Big|\frac{1}{2\tau}\langle(I-\S_L^2(\tau/2))u^n, u^n\rangle-\frac{1}{2}\langle\epsilon^2 (-\Delta)u^n, u^n\rangle\Big|\\
\le&\Big|\frac{1}{2\tau}\langle(I-\S_L^2(\tau/2))u^n, u^n\rangle-\frac{1}{2\tau}\langle(I-e^{\epsilon^2\tau L})u^n, u^n\rangle\Big|+\Big|\frac{1}{2\tau}\langle(1-e^{\epsilon^2\tau L})u^n, u^n\rangle-\frac{1}{2}\langle\epsilon^2 (-L)u^n, u^n\rangle\Big|\\
&+\Big|\frac{1}{2}\langle\epsilon^2 (-L)u^n, u^n\rangle-\frac{1}{2}\langle\epsilon^2 (-\Delta)u^n, u^n\rangle\Big|\\
\lesssim&\O(\tau).
\end{align*}
This completes our proof.
\end{proof}

\subsubsection{Nonlinear-Linear-Nonlinear Strang splitting}
\label{subsec:nln_splitting}
We now consider the more practical nonlinear-linear-nonlinear splitting $u^{n+1}=\S_N(\tau/2)\S_L(\tau)\S_N(\tau/2)u^n$. If the energy in the previous section is directly used, the proof fails. Therefore, a new non-negative energy is needed, which should ideally converge to the exact Ginzburg--Landau energy. Following the logical arguments in~\cite{StrangSplittingAC} we define a suitable $\Tilde{u}$ as follows in the proof below. \delete{Since the proof in~\cite{StrangSplittingAC} uses exact spatial and temporal update for the linear part, minor additional requirements are needed in order for the proof to hold. Notice that compared to the original proof in~\cite{StrangSplittingAC}, we further assumed that the spatial discretization size is fixed in order to conduct our analysis.}
\begin{theorem}[Modified Energy Dissipation, nonlinear-linear-nonlinear]
\label{thm:nln_energy}
For any symmetric positive-definite $S_L(\tau)$ with spectrum in $(0, 1]$, define the modified energy as
\begin{equation}
\label{eqn:nln_energy}
\Tilde{E}(u^n):=\frac{1}{2\tau}\Big\langle(I-S_L(\tau))S_N(\tau/2)u^n,S_L(\tau)S_N(\tau/2)u^n\Big\rangle+\langle g(S_L(\tau)S_N(\tau/2)u^n),1\rangle,
\end{equation}
where the function $g$ is defined as in~\Cref{eqn:modified_nonlinear_energy}. Then under the nonlinear-linear-nonlinear Strang splitting defined above, we have
\begin{enumerate}
\item $\Tilde{E}(u^{n+1})\le\Tilde{E}(u^{n})$,
\item For fixed grid size $h$, $\Tilde{E}(u^n)$ converges to $E(u^n)$ uniformly as $\tau\to0$; more precisely, for $0\le\tau\le\min\{1, \text{CFL constraint}\}$, assuming enough regularity of $u^0$ we have that
$$\underset{n\ge0}{\sup}|\Tilde{E}(u^n)-E(u^n)|\lesssim\tau.$$
\end{enumerate}
\end{theorem}
\begin{remark}
As will be shown in~\Cref{subsubsec:modified_rkl_rkg} and~\Cref{subsubsec:positive_definite_rkl_rkg}, the assumptions can be satisfied by modified versions of RKL1 and RKG1. It will also be naturally satisfied by RKL2, but requires some more careful estimations and restrictions for RKG2. 
\end{remark}
\begin{proof}
The non-negativity of $\Tilde{E}(u^n)$ can be easily deduced from the non-negativity of $g$ and the fact that $(I-S_L(\tau))S_L(\tau)$ is positive semi-definite. To show that the modified energy dissipates, set $\Tilde{u}^n:=\S_L(\tau)\S_N\Big(\frac{\tau}{2}\Big)u^n$. We know that
$$\Tilde{u}^{n+1}=\S_L(\tau)\S_N(\tau)\Tilde{u}^n,$$
so that
$$\S_L(\tau)^{-1}\Tilde{u}^{n+1}=\frac{e^{\tau}\Tilde{u}^n}{\sqrt{1+(e^{2\tau}-1)(\Tilde{u}^n)^2}}.$$
Rewriting the above and using direct computation, we find
$$\frac{1}{\tau}(\S_L(\tau)^{-1}-I)\Tilde{u}^{n+1}+\frac{1}{\tau}(\Tilde{u}^{n+1}-\Tilde{u}^n)=-g'(\Tilde{u}^n).$$
Also by direct computation and self-adjointness of $\S_L(\tau)$ and the fact that all the eigenvalues of $\S_L(\tau)$ are non-negative, we have that
\begin{align*}
    &\frac{1}{\tau}\langle(\S_L(\tau)^{-1}-I)\Tilde{u}^{n+1}, \Tilde{u}^{n+1}-\Tilde{u}^n\rangle\\
    =&\frac{1}{\tau}\langle(I-\S_L(\tau))\S_N\Big(\frac{\tau}{2}\Big)u^{n+1}, \S_L(\tau)\S_N\Big(\frac{\tau}{2}\Big)(u^{n+1}-u^n)\rangle\\
    =&\frac{1}{\tau}\langle(I-\S_L(\tau))\sqrt{\S_L(\tau)}\S_N\Big(\frac{\tau}{2}\Big)u^{n+1}, \sqrt{\S_L(\tau)}\S_N\Big(\frac{\tau}{2}\Big)(u^{n+1}-u^n)\rangle\\
    =&\frac{1}{2\tau}\Big(\langle(I-\S_L(\tau))\sqrt{\S_L(\tau)}\S_N\Big(\frac{\tau}{2}\Big)u^{n+1}, \sqrt{\S_L(\tau)}\S_N\Big(\frac{\tau}{2}\Big)u^{n+1}\rangle\\
    &-\langle(I-\S_L(\tau))\sqrt{\S_L(\tau)}\S_N\Big(\frac{\tau}{2}\Big)u^n, \sqrt{\S_L(\tau)}\S_N\Big(\frac{\tau}{2}\Big)u^n\rangle\\
    &+\langle(I-\S_L(\tau))\sqrt{\S_L(\tau)}\S_N\Big(\frac{\tau}{2}\Big)(u^{n+1}-u^n), \sqrt{\S_L(\tau)}\S_N\Big(\frac{\tau}{2}\Big)(u^{n+1}-u^n)\rangle\Big)\\
    =&\frac{1}{2\tau}\Big(\langle(\S_L(\tau)^{-1}-I)\Tilde{u}^{n+1},\Tilde{u}^{n+1}\rangle-\langle(\S_L(\tau)^{-1}-I)\Tilde{u}^{n}, \Tilde{u}^n\rangle\\
    &+\langle(I-\S_L(\tau))\sqrt{\S_L(\tau)}\S_N\Big(\frac{\tau}{2}\Big)(u^{n+1}-u^n), \sqrt{\S_L(\tau)}\S_N\Big(\frac{\tau}{2}\Big)(u^{n+1}-u^n)\rangle\Big).
\end{align*}
Multiplying the equation and summing over the grid, and using Taylor expansion of $g$ up to second-order, we have that
\begin{align*}
&\Tilde{E}(\Tilde{u}^{n+1})-\Tilde{E}(\Tilde{u}^{n})\\
=&-\frac{1}{2\tau}\langle(I-\S_L(\tau))\sqrt{\S_L(\tau)}\S_N\Big(\frac{\tau}{2}\Big)(u^{n+1}-u^n), \sqrt{\S_L(\tau)}\S_N\Big(\frac{\tau}{2}\Big)(u^{n+1}-u^n)\rangle\\
&-\Big\langle\Big(\frac{1}{\tau}-\frac{1}{2}g''(\xi^n)\Big)(\Tilde{u}^{n+1}-\Tilde{u}^n)^2, 1\Big\rangle.
\end{align*}
By noticing that $I-\S_L(\tau)$ is symmetric positive-definite and noticing that $\frac{1}{\tau}-\frac{1}{2}g''(\xi^n)\ge0$ \delete{(detailed derivation in the paper~\cite{StrangSplittingAC})}, we can conclude that $\Tilde{E}(\Tilde{u}^{n+1})\le\Tilde{E}(\Tilde{u}^{n})$ as desired.

Now we show that the modified energy converges to the exact energy. We will first show that $\Tilde{u}^n$ and $u^n$ differ by $\mathcal{O}(\tau)$. By definition we have that $\Tilde{u}^n:=\S_L(\tau)\S_N\Big(\frac{\tau}{2}\Big)u^n$; taking the difference and using Taylor expansion, and by noticing that the norm of discrete Laplacian $\|L\|_\infty$ is bounded for fixed grid size, we have that
\begin{align*}
\|\Tilde{u}^n-u^n\|_\infty=&\|\S_L(\tau)\S_N\Big(\frac{\tau}{2}\Big)u^n-u^n\|_\infty\\
=&\|(I+\tau L+\mathcal{O}(\tau ^2))(e^{\tau/2}+\mathcal{O}(\tau))u^n-u^n\|_\infty\\
=&\|(I+\tau L+\mathcal{O}(\tau ^2))(I+\mathcal{O}(\tau))u^n-u^n\|_\infty\\  
\lesssim&\,\mathcal{O}(\tau).
\end{align*}
Now notice that $\|\Tilde{u}^n\|_\infty\le\|S_L(\tau/2)\|_\infty\|S_N(\tau/2)u^n\|_\infty\le C\|S_L(\tau/2)\|_2\le C$ for some grid size dependent constant $C$ uniformly when the CFL condition is satisfied. Then by the same argument as in equations (2.26) and (2.27) of~\cite{StrangSplittingAC}, we are done with the nonlinear part convergence.
\\
For the inner product part of the energy, by direct calculation we know that
\begin{align*}
&\left|\frac{1}{2\tau}\langle(I-\S_L(\tau))\sqrt{\S_L(\tau)}\S_N\Big(\frac{\tau}{2}\Big)u^n, \sqrt{\S_L(\tau)}\S_N\Big(\frac{\tau}{2}\Big)u^n\rangle-\frac{1}{2}\langle\epsilon^2 Lu^n, u^n\rangle\right|\\
\le&\Big|\frac{1}{2\tau}\langle(I-\S_L(\tau))\sqrt{\S_L(\tau)}\S_N\Big(\frac{\tau}{2}\Big)u^n, \sqrt{\S_L(\tau)}\S_N\Big(\frac{\tau}{2}\Big)u^n\rangle-\frac{1}{2}\langle\epsilon^2 (-L)\sqrt{\S_L(\tau)}\S_N\Big(\frac{\tau}{2}\Big)u^n, \sqrt{\S_L(\tau)}\S_N\Big(\frac{\tau}{2}\Big)u^n\rangle\Big|\\
&+\Big|\frac{1}{2}\langle\epsilon^2(-L)\sqrt{\S_L(\tau)}\S_N\Big(\frac{\tau}{2}\Big)u^n, \sqrt{\S_L(\tau)}\S_N\Big(\frac{\tau}{2}\Big)u^n\rangle-\frac{1}{2}\langle\epsilon^2(-L)u^n, u^n\rangle \Big|\\
&+\Big|\frac{1}{2}\langle\epsilon^2(-L)u^n, u^n\rangle-\frac{1}{2}\langle\epsilon^2(-\Delta)u^n, u^n\rangle\Big|\lesssim\O(\tau). 
\end{align*}
This completes our proof.
\end{proof}
\begin{remark}
The linear-nonlinear-linear energy dissipation utilizes the power of Strang splitting instead of the power of super-time-stepping methods. In fact, any linear stable approximation to the heat equation can be fitted into the proof outline above. The nonlinear-linear-nonlinear Strang splitting energy dissipation, on the other hand, does rely on the fact that the stability polynomial of the second-order SRK methods are (almost) positive-definite.
\end{remark}
It is easy to verify all the conditions required in the theorem above, except for the positive-definiteness of the $S_L(\tau)$ operators for the RKL and RKG methods. We therefore proceed and check the condition in the following two sections, and show that the RKL and RKG methods are indeed able to be fitted into the proof above up to modifications and minor additional timestep constraint.

\subsubsection{Positive-definite modified RKL1 and RKG1 methods}
\label{subsubsec:modified_rkl_rkg}
Since the nonlinear-linear-nonlinear energy definition and the proof rely on the fact that $S_L(\tau)$ has spectrum within $(0, 1]$, we now  determine the condition that such an assumption holds. For the RKL1 and RKG1 methods
\begin{align*}
u^{t+\tau}=P_{s}\Big(I+\frac{2\tau}{s^2+s}\bm{M}\Big)u^t,
\end{align*}
\begin{align*}
u^{t+\tau}=\frac{2}{(s+1)(s+2)}C^{(\lambda)}_{s}\Big(I+\frac{4\tau}{s^2+3s}\bm{M}\Big)u^t.
\end{align*}
According to the Legendre and Gegenbauer polynomial properties, $|S_L(\tau)|\le1$ is the best possible estimation. Therefore, modifications are needed in order for the first-order methods to have a non-negative spectrum. We then follow the general derivation procedure as in~\Cref{subsec:rkl_rkg} to derive new modified first-order RKL and RKG methods.

For the modified RKL1 method, instead of using the Legendre polynomial, consider the shifted Legendre polynomial $\Tilde{P}_n(x)=\frac{1+P_n(x)}{2}$, where $P_n(x)$ is the $n$-th order Legendre polynomial. Since the Legendre polynomials satisfy $P_n(x)\in[-1, 1]$ for $x\in[-1, 1]$ and $\forall n\in\N$, we have that $\Tilde{P}_n(x)\in[0, 1]$ for $x\in[-1, 1]$ and $\forall n\in\N$. The three-point recurrence relationship is then
$$j\,\Tilde{P}_j(x)=(2j-1)\,x\,\Tilde{P}_{j-1}(x)-(j-1)\,\Tilde{P}_{j-2}(x)+\frac{2j-1}{2}\,(1-x),\qquad j\ge 2.$$
Now assume that the stability polynomial takes the form of
$$R_s(z)=a_s+b_s\Tilde{P}_n(w_0+w_1z).$$
For first-order, we can assume that $a_s=0$ for all $s\in\N$ and choose $w_0=1$ as in~\cite{RKL2}, so that $$R_s(z)=\Tilde{P}_s\Big(1+\frac{4}{s^2+s}z\Big).$$
Plugging the stability polynomial into the three-point recurrence relationship, we obtain that
$$\Tilde{P}_j\Big(1+\frac{4}{s^2+s}z\Big)=\frac{2j-1}{j}\Big(1+\frac{4}{s^2+s}z\Big)\Tilde{P}_{j-1}\Big(1+\frac{4}{s^2+s}z\Big)-\frac{j-1}{j}\Tilde{P}_{j-2}(x)+\frac{2j-1}{2j}\Big(-\frac{4}{s^2+s}z\Big).$$
The modified RKL1 method is then
\begin{align*}
\begin{split}
Y_0&=u(t_0),\\
Y_1&=Y_0+\Tilde{\mu}_1\tau\bm{M}Y_0,\\
Y_j&=\mu_jY_{j-1}+\nu_jY_{j-2}+\Tilde{\mu}_j\tau\bm{M}Y_{j-1}+\Tilde{\gamma}_j\tau\bm{M}Y_0,\\
u(t_0+\tau)&=Y_s,
\end{split}
\end{align*}
where
\begin{align*}
\mu_j=\frac{2j-1}{j}, \nu_j=\frac{1-j}{j}, \Tilde{\mu}_j=\frac{2j-1}{j}\frac{4}{s^2+s}, \Tilde{\gamma}_j=-\frac{4j-2}{js(s+1)}.
\end{align*}

Similarly, we can derive the modified RKG1 method by using the shifted Gegenbauer polynomial 
$$\Tilde{C}^{(\lambda)}_n(x)=\frac{1}{2}\Big(1+\frac{1}{(n+1)(n+2)}C^{(\lambda)}_n(x)\Big).$$
By using exactly the same procedure as above, the RKG1 method can be written as
\begin{align*}
\begin{split}
Y_0&=u(t_0),\\
Y_1&=Y_0+\Tilde{\mu}_1\tau\bm{M}Y_0,\\
Y_j&=\mu_jY_{j-1}+\nu_jY_{j-2}+\Tilde{\mu}_j\tau\bm{M}Y_{j-1}+\Tilde{\gamma}_j\tau\bm{M}Y_0,\\
u(t_0+\tau)&=Y_s,
\end{split}
\end{align*}
where
\begin{align*}
\begin{split}
\mu_j=\frac{2j+1}{j+2}, \nu_j=-\frac{j-1}{j+2}, \Tilde{\mu}_j=\frac{2j+1}{j+2}\frac{8}{s^2+3s}, \Tilde{\gamma}_j=-\frac{2j+1}{(2j+2)}\frac{8}{s^2+3s},\Tilde{\mu}_1=\frac{4}{s^2+3s}.
\end{split}
\end{align*}

The advantage of these modified RKL1 and RKG1 is that they control the eigenvalues of the stability polynomial to be all non-negative, so that they can be directly fitted into the proof of~\Cref{thm:nln_energy}. Moreover, those methods naturally satisfy the discrete maximum principle since the regular RKL1 and RKG1 methods are, under the corresponding boundary condition requirements of the RKL and RKG methods. The drawback of this modified RKL1 is that the maximum possible $\tau$ is only half of the original RKL1 by following the same derivation as in~\cite{RKL2}, although it still follows the $\O(s^2)$ order.

\subsubsection{Positive-definiteness estimation for RKL2 and RKG2 methods}
\label{subsubsec:positive_definite_rkl_rkg}
We now turn our attention to the second-order RKL and RKG. According to the derivations in~\cite{RKL2} and~\cite{SKARAS_RKG2021109879}, it is not difficult to see that the shifted RKL2 and RKG2 are the same as the usual RKL2 and RKG2, since there are exactly three equations and three unknowns. Therefore, we need to estimate the conditions for RKL2 and RKG2 to give us positive (semi)-definite stability polynomials. 

Recall that the stability polynomial for RKL2 is 
$$R_s(\tau\bm{M})=\Big(1-\frac{s^2+s-2}{2s(s+1)}\Big)I+\frac{s^2+s-2}{2s(s+1)}P_{s}\Big(I+\frac{4\tau}{s^2+s-2}\bm{M}\Big).$$
Notice that $0\le\frac{s^2+s-2}{2s(s+1)}<\frac{1}{2}$ for $s\in\Z_+$; since the Legendre polynomial is bounded by 1 when the argument is within $[-1, 1]$, we know that $0<R_s(\tau\bm{M})\le1$ when $\tau$ is within the stability constraint. This shows that RKL2 automatically gives us a positive-definite stability polynomial.

The RKG2 method does not automatically have a positive-definite stability polynomial. Recall that the RKG2 stability polynomial is
$$\Big(1-\frac{2(s-1)(s+4)}{3s(s+3)}\Big)I+\frac{4(s-1)(s+4)}{3s(s+1)(s+2)(s+3)}C^{(\lambda)}_{s}\Big(I+\frac{6\tau}{(s+4)(s-1)}\bm{M}\Big).$$
Notice that $0\le\frac{2(s-1)(s+4)}{3s(s+3)}<\frac{2}{3}$ for $s\in\Z_+$; since the Gegenbauer polynomial of degree $s$ and $\lambda=\frac{3}{2}$ has value within $\Big[-\frac{(s+1)(s+2)}{2}, \frac{(s+1)(s+2)}{2}\Big]$, the spectrum of the RKG2 stability polynomial takes values within $(-\frac{1}{3}, 1)$. Since the most negative values of the Gegenbauer polynomials correspond to argument near $-1$ with $s$ odd, we need to estimate the domain of the Gegenbauer polynomial so that the argument within the polynomial is away from $-1$. With more asymptotic analysis and the control of the positions for the zeros of the Gegenbauer polynomial, we manage to show that the characteristic polynomial of RKG2 is \textit{``almost positive-definite''}\delete{, with the precise result formulated in the following theorem.}\revision{ in the sense that with a small sacrifice of the discrete timestep size away from the maximum possible discrete timestep size, the RKG2 operator is positive-definite. More precisely, it is shown in the theorem below.}
\begin{theorem}[RKG2 operator is ``almost positive-definite'']
Consider RKG2 in~\Cref{eqn:RKG2} with parameter $\lambda=\frac{3}{2}$, and let $\tau_{\max}$ denote the largest timestep satisfying the CFL constraint. Then for any $0<\epsilon\ll1$, there exists $s_0\in\N$ such that for all $0<\tau\le(1-\epsilon)\tau_{\max}$ and $s\ge s_0$, $S_L(\tau)$ is positive-definite.
\end{theorem}
\begin{proof}
First recall that $C^{(\lambda)}_s(x)=\frac{\Gamma(\lambda+\frac{1}{2})}{\Gamma(2\lambda)}\frac{\Gamma(s+2\lambda)}{\Gamma(s+\lambda+\frac{1}{2})}P_s^{(\lambda-\frac{1}{2}, \lambda-\frac{1}{2})}(x)$ where $\Gamma(x)$ is the standard gamma function~\cite{szego1975orthogonal}. Now arbitrarily fix $\delta>0$. Now fix $c>0$ small enough such that $\cos\frac{c}{s}>x_s$ where $x_s$ is the right-most root of $C^{(\lambda)}_s(x)=C^{(3/2)}_s(x)$. Notice that this is attainable via Theorem 4.1 of~\cite{gatteschi1985zeros}: for $\vartheta_{s, k}(1, 1)$ such that $\cos(\vartheta_{s, k}(1, 1))$ are the zeros of the Gegenbauer polynomial $P^{(1, 1)}_s(x)=C^{(3/2)}_s(x)$, we have that specifically
$$\vartheta_{s, 1}(1, 1)=\frac{j_{1, 1}}{\nu}+\O(s^{-5}),$$
where $j_{1, 1}$ denotes the left-most zero of $J_1(x)$ for $x\in[0, \infty)$, and $\nu=\left((s+\frac{3}{2})^2-\frac{1}{4}\right)^{1/2}\sim O(s)$; since $\cos(\vartheta_{s, 1}(1, 1))$ corresponds to the right-most zero of $C^{(3/2)}_s(x)$, we can choose $0<c\lesssim j_{1, 1}$ small enough to achieve the desired goal.

Notice that the positive-definiteness of the RKG2 operator, with the timestep restriction in the theorem statement, is equivalent to showing that
\begin{equation}
\label{eqn:normalized_gegenbauer_condition}
\Tilde{C}^{(\lambda)}_s(x):=\frac{C^{(\lambda)}_s(x)}{C^{(\lambda)}_s(1)}>-\frac{1}{2},
\end{equation}
for $x\in[-1, 1]$ that is not close to $-1$, and as $s\to\infty$. This is because the eigenvalues lie in the interval of the form
$$1-\frac{2(s-1)(s+4)}{3s(s+3)}+\frac{4(s-1)(s+4)}{3s(s+1)(s+2)(s+3)}C^{(\lambda)}_s(x)=1-\frac{2(s-1)(s+4)}{3s(s+3)}+\frac{2(s-1)(s+4)}{3s(s+3)}\Tilde{C}^{(\lambda)}_s(x)\delete{.}\revision{,}$$
along with the fact that $\frac{2(s-1)(s+4)}{3s(s+3)}\underset{s\to\infty}{\longrightarrow}\frac{2}{3}$. Therefore, we will show that~\cref{eqn:normalized_gegenbauer_condition} holds; more specifically, we will show that after the standard change of variable $x:=\cos\theta$, it holds for $0<\theta\le\pi-\delta$ in order to complete our proof.

To conduct the estimation, we dissect the $\theta$ region $[0, \pi-\delta]$ into three different regimes. We first choose $d>0$ to be specified later.
\paragraph{Regime 1: $0\le \theta\le\frac{c}{s}$} For the region for $0<\theta<\frac{c}{s}$, by the way we construct $c$, naturally $C^{(\lambda)}_s(x)\ge0$ so that $\Tilde{C}^{(\lambda)}_s(x)\ge0$ since $C^{(\lambda)}_s(1)>0$ always holds.

\paragraph{Regime 2: $\frac{c}{s}\le\theta\le\frac{d}{s}$} In this regime, consider the Mehler–Heine formula~\cite{szego1975orthogonal}
$$\underset{s\to\infty}{\lim}s^{-\alpha}P^{(\alpha, \beta)}\left(\cos\left(\frac{z}{s}\right)\right)=\left(\frac{z}{2}\right)^{-\alpha}J_\alpha(z).$$
Notice that
$$\Tilde{C}^{(3/2)}_s(x)=\frac{C^{(3/2)}_s(x)}{C^{(3/2)}_s(1)}=\frac{s+2}{2}P_s^{(1, 1)}(x)\frac{2}{(s+1)(s+2)}=\frac{1}{s+1}P_s^{(1, 1)}(x).$$
Therefore, we have that 
$$\underset{s\to\infty}{\lim}\Tilde{C}^{(3/2)}_s\left(\cos\left(\frac{z}{s}\right)\right)=\frac{2J_1(z)}{z}.$$
By elementary estimation of $\frac{J_1(z)}{z}$, we have that for all sufficiently large $s>0$, $\Tilde{C}^{(3/2)}_s\left(\cos\left(\frac{z}{s}\right)\right)>-\frac{1}{2}$ for all $z\in[c, d]$, or equivalently $\theta\in[\frac{c}{s}, \frac{d}{s}]$.
\paragraph{Regime 3: $\frac{d}{s}\le\theta\le\pi-\delta$} By theorem 8.21.12 in~\cite{szego1975orthogonal},
\begin{align*}
\begin{split}
&\left(\sin\frac{\theta}{2}\right)^{\alpha}
\left(\cos\frac{\theta}{2}\right)^{\beta}
P_s^{(\alpha,\beta)}(\cos\theta)=N^{-\alpha}\frac{\Gamma(s+\alpha+1)}{s!}\Big(\frac{\theta}{\sin\theta}\Big)^{1/2}J_\alpha(N\theta)+\theta^{1/2}\O(s^{-3/2}),\quad\frac{c}{s}\le\theta\le\pi-\delta,
\end{split}
\end{align*}
where $N=s+\frac{\alpha+\beta+1}{2}$, $J_\alpha(x)$ is the Bessel function of the first kind of order $\alpha$~\cite{szego1975orthogonal}.
In our case of $\lambda=\frac{3}{2}$, $\alpha=\lambda-\frac{1}{2}=1$, and the above formula is equivalent to
\begin{align*}
\begin{split}
&\frac{1}{2}\sin\theta\,P_s^{(1,1)}(\cos\theta)=N^{-1}\frac{\Gamma(s+2)}{s!}\Big(\frac{\theta}{\sin\theta}\Big)^{1/2}J_{1}(N\theta)+\theta^{1/2}\O(s^{-3/2}),\quad\frac{d}{s}\le\theta\le\pi-\delta.
\end{split}
\end{align*}
Using the relationship between $\Tilde{C}^{(3/2)}_s(x)$ and $P^{(1, 1)}(x)$, we have that
$$\Tilde{C}^{(3/2)}_s(x)=\frac{2}{\sin\theta(s+1)}N^{-1}\frac{\Gamma(s+2)}{s!}\Big(\frac{\theta}{\sin\theta}\Big)^{1/2}J_{1}(N\theta)+\theta^{1/2}\O(s^{-3/2}),\quad\frac{d}{s}\le\theta\le\pi-\delta.$$
By using the fact that $|J_1(z)|\le\frac{C}{z^{1/2}}$ for some $C>0$, we have that
$$\left|\frac{2}{\sin\theta(s+1)}N^{-1}\frac{\Gamma(s+2)}{s!}\Big(\frac{\theta}{\sin\theta}\Big)^{1/2}J_{1}(N\theta)\right|\le A\frac{1}{\sin\theta(s+1)}\Big(\frac{\theta}{\sin\theta}\Big)^{1/2}\frac{1}{(s\theta)^{1/2}}\le\frac{A}{s^{3/2}(\sin\theta)^{3/2}},$$
for some constant $A$. 
Now choose $d>c$ such that $\frac{A}{d^{3/2}}<\frac{1}{4}$; notice this bound is not circular, as $A$ is independent of $d$ but depends only on the Bessel function bound uniform constant $C$. The worst case scenario of the bound above occurs at $\frac{d}{s}$ for all $s\ge s_0$ with some $s_0\in\N$ large enough, so that $$\left|\frac{2}{\sin\theta(s+1)}N^{-1}\frac{\Gamma(s+2)}{s!}\Big(\frac{\theta}{\sin\theta}\Big)^{1/2}J_{1}(N\theta)\right|\le \frac{A}{d^{3/2}}.$$
Then for $s\gg1$ such that the error term is also less than $\frac{1}{4}$ in magnitude and the leading term dominates, we get that
$$|\Tilde{C}^{(3/2)}_s(x)|<\frac{1}{2}.$$
This ends our three regime analysis.

Then since we need $\theta\in[0, \pi-\delta]$, we have that $x=\cos\theta\in[-\cos\delta, 1]$. By the RKG2 argument definition, we have that 
$$x(\tau)_j=1+\frac{6\tau}{(s+4)(s-1)}{\mu}_j,$$
where $\mu_j$ is an eigenvalue of the discrete Laplacian operator. Since we want $x\ge-\cos\delta$, we need that
$$1+\frac{6\tau}{(s+4)(s-1)}{\mu}_{\max}\ge-\cos\delta\Longrightarrow\tau\le\frac{(s+4)(s-1)}{6\mu_{\max}}(1+\cos\delta),$$
where $\mu_{\max}$ is the maximum eigenvalue in absolute value. Notice that the $\tau_{\max}$ of the RKG2 method is achieved at $\tau=\frac{(s+4)(s-1)}{3\mu_{\max}}$ because we want to set $x(\tau)=-1$ in this case; therefore, we get that
$$\tau\le\frac{1+\cos\delta}{2}\tau_{\max}=(1-\epsilon)\tau_{\max}\revision{,}$$
if we set $\epsilon:=\frac{1-\cos\delta}{2}>0$. As desired, since $\delta>0$ is a small number, we have that $\epsilon>0$ is a small number. This concludes our proof.
\end{proof}

\begin{remark}
The above theorem is only a rough estimate; nevertheless, it shows that RKG2 can be controlled away from its most negative value, so that RKG2 with a more restricted timestepping and large enough stage number $s$ can be fitted into the proof of~\Cref{thm:nln_energy}.
\end{remark}

\section{Experimental results}
\label{sec:experiment}

\begin{table}[!t]
    \begin{center}
    \begin{tabular}{|c|c|c|c|c|c|c|c|c|}
    \hline
    \multirow{3}{*}{$\tau=$0.25}&\multicolumn{4}{|c|}{RKL2}&\multicolumn{4}{|c|}{RKG2}\\
    \cline{2-9}
    &\multicolumn{2}{|c|}{$L^2$-norm}&\multicolumn{2}{|c|}{$L^\infty$-norm}&\multicolumn{2}{|c|}{$L^2$-norm}&\multicolumn{2}{|c|}{$L^\infty$-norm}\\
    \cline{2-9}
    & Error & Rate & Error & Rate & Error & Rate & Error & Rate\\
    \hline
    $\tau$ & - & - & - & - & - & - & - & -\\
    \hline
    $\tau/2$ & 1.312e-05 & - & 9.263e-08 & - & 1.304e-05 & - & 9.215e-08 & - \\
    \hline
    $\tau/4$ & 3.301e-06 & 1.991 & 2.331e-08 & 1.991 & 3.280e-06 & 1.991 & 2.319e-08 & 1.990 \\ 
    \hline
    $\tau/8$ & 8.264e-07 & 1.998 & 5.837e-09 & 1.998 & 8.214e-07 & 1.998 & 5.808e-09 & 1.998 \\
    \hline
    $\tau/16$ & 2.067e-07 & 1.999 & 1.460e-09 & 1.999 & 2.054e-07 & 1.999 & 1.452e-09 & 1.999 \\
    \hline
    $\tau/32$ & 5.168e-08 & 2.000 & 3.650e-10 & 2.000 & 5.136e-08 & 2.000 & 3.632e-10 & 2.000 \\
    \hline
    \end{tabular}
    \end{center}
    \vspace{-1em}
    \caption{Temporal rate of convergence estimation for the RKL2 and RKG2 methods applied to the Allen--Cahn equation, along with Strang splitting. The convergence rates for both methods are close to 2 as expected.}
    \label{table:order_estimation_allen_cahn}
\end{table}

\begin{table}[!t]
    \begin{center}
    \begin{tabular}{|c|c|c|c|c|c|c|c|c|}
    \hline
    \multirow{3}{*}{$\tau=$1.6e-3}&\multicolumn{4}{|c|}{RKL2}&\multicolumn{4}{|c|}{RKG2}\\
    \cline{2-9}
    &\multicolumn{2}{|c|}{$L^2$-norm}&\multicolumn{2}{|c|}{$L^\infty$-norm}&\multicolumn{2}{|c|}{$L^2$-norm}&\multicolumn{2}{|c|}{$L^\infty$-norm}\\
    \cline{2-9}
    & Error & Rate & Error & Rate & Error & Rate & Error & Rate\\
    \hline
    $\tau$ & - & - & - & - & - & - & - & -\\
    \hline
    $\tau/2$ & 8.529e-08 & - & 5.288e-10 & - & 6.404e-08 & - & 3.971e-10 & - \\
    \hline
    $\tau/4$ & 2.135e-08 & 1.998 & 1.324e-10 & 1.998 & 1.603e-08 & 1.999 & 9.942e-11 & 1.998 \\ 
    \hline
    $\tau/8$ & 5.340e-09 & 1.999 & 3.313e-11 & 1.999 & 4.009e-09 & 1.999 & 2.488e-11 & 1.999 \\
    \hline
    $\tau/16$ & 1.335e-09 & 2.000 & 8.289e-12 & 1.999 & 1.002e-09 & 2.000 & 6.222e-12 & 1.999 \\
    \hline
    $\tau/32$ & 3.339e-10 & 2.000 & 2.080e-12 & 1.995 & 2.504e-10 & 2.001 & 1.558e-12 & 1.998 \\
    \hline
    \end{tabular}
    \end{center}
    \vspace{-1em}
    \caption{Temporal rate of convergence estimation for the RKL2 and RKG2 method applied to the Cahn--Hilliard equation. The convergence rates for both methods are close to 2 as expected.}
    \label{table:order_estimation_cahn_hilliard}
\end{table}

In this section, we present numerical experiments of the RKL2 and RKG2 methods applied to the Allen--Cahn and Cahn--Hilliard equations. \delete{We assume that all the Allen--Cahn Strang splitting takes the order of nonlinear-linear-nonlinear splitting as in~\mbox{\Cref{subsec:nln_splitting}} since the RKL and RKG methods are only required once per step.}\revision{For all the Allen--Cahn simulations, we take the nonlinear-linear-nonlinear Strang splitting order written in}~\Cref{subsec:nln_splitting}. We first verify the second-order convergence rate of the RKL2 and RKG2 methods when applied to the Allen--Cahn and Cahn--Hilliard equations. Following that, we present qualitative results for the Allen--Cahn and Cahn--Hilliard coarsening behaviors, by comparing them with known results in existing papers and random initial condition coarsening. Furthermore, we simulate the Allen--Cahn and the Cahn--Hilliard equations for an $L$-shaped domain\delete{to demonstrate its superiority over other Fourier spectral method based time marching schemes.}\revision{; finite difference approximation rather than the Fourier spectral method is used for spatial discretization.}

\subsection{Convergence rate test}
We consider the two dimensional Allen--Cahn equation in~\Cref{eqn:AllenCahn}, with $\epsilon=0.1$. The simulation domain is $[0, 2\pi]^2$ with periodic boundary condition. The spatial domain is divided into $256\times256$ intervals, so the spatial discretization error is smaller than  the time discretization error. The initial condition is 
$$u(x, y, 0)=0.05\sin(x)\sin(y).$$
We test the order of convergence for the second-order Runge--Kutta-Legendre (RKL2) and Runge--Kutta--Gegenbauer (RKG2) methods, when using Strang splitting. We choose the number of substeps $s=100$ and simulate until $T=1$. We adopt the standard step-doubling convergence rate estimation since there is no exact solution available; 
$$\text{convergence rate}\approx\frac{\log(\|e_i\|/\|e_{i+1}\|)}{\log2},$$
where $e_i:=u_{\tau/2^i}(T)-u_{\tau/2^{i+1}}(T)$ if we denote $u_\tau(T)$ to be the simulated solution with step size $\tau$ at time $T$, and $\|\cdot\|$ is either $L^2$-norm or $L^\infty$-norm. As shown in~\Cref{table:order_estimation_allen_cahn}, the rate of convergence for the RKL2 and RKG2 methods in both norms are close to the expected second-order.
We test the order of convergence for RKL2 and RKG2 when applied to the Cahn--Hilliard equation, under the same conditions as above up to $T=0.016$. The temporal convergence rate is shown in~\Cref{table:order_estimation_cahn_hilliard}, with the desired second-order convergence.

\begin{figure}[!t]
    \centering
    \begin{minipage}{0.8\linewidth}
        \begin{subfigure}{\linewidth}
        \captionsetup{justification=centering}
        \includegraphics[width=0.2\linewidth]{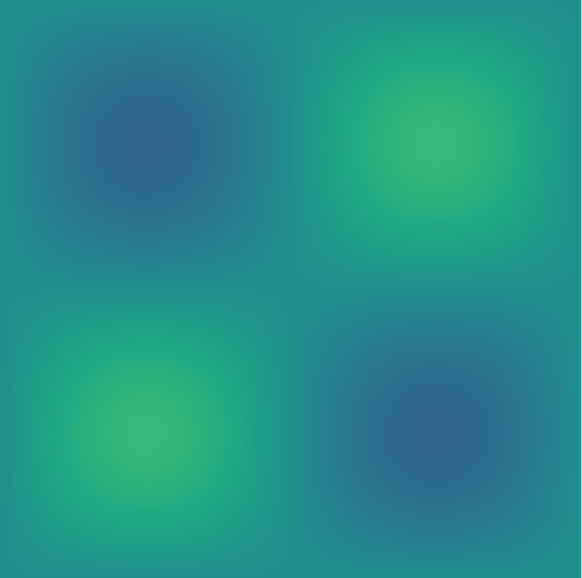}
        \hspace{0.7cm}
        \includegraphics[width=0.2\linewidth]{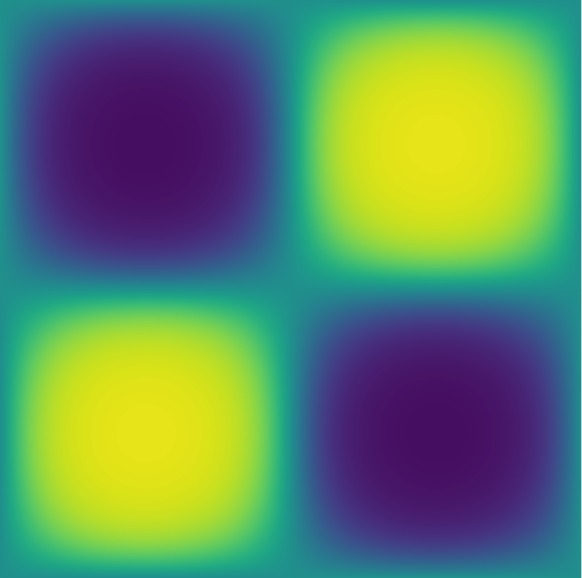}
        \hspace{0.7cm}
        \includegraphics[width=0.2\linewidth]{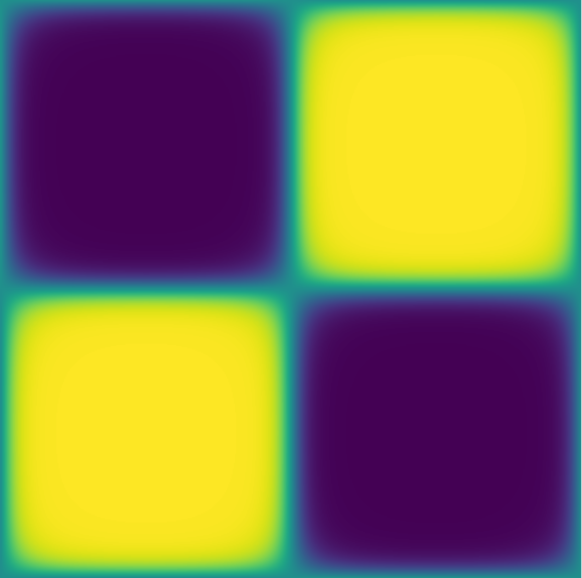}
        \hspace{0.7cm}
        \includegraphics[width=0.2\linewidth]{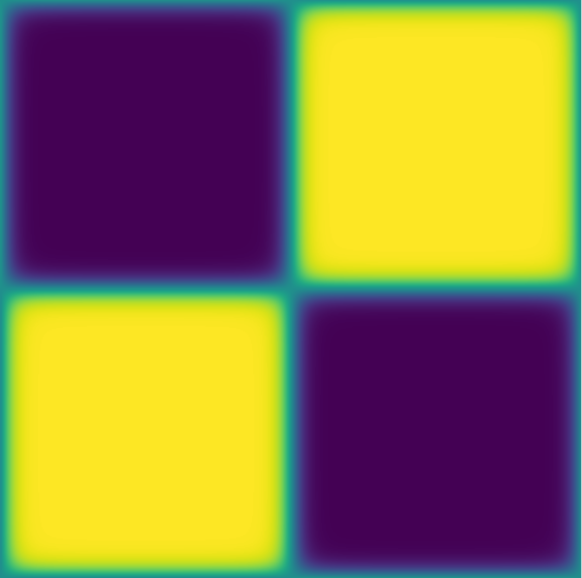}
        \caption{Allen--Cahn simulation, RKL2}
        \end{subfigure}
        \begin{subfigure}{\linewidth}
        \captionsetup{justification=centering}
        \includegraphics[width=0.2\linewidth]{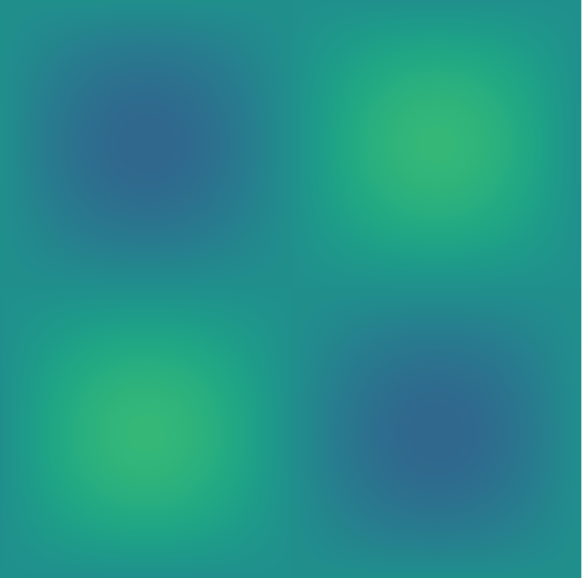}
        \hspace{0.7cm}
        \includegraphics[width=0.2\linewidth]{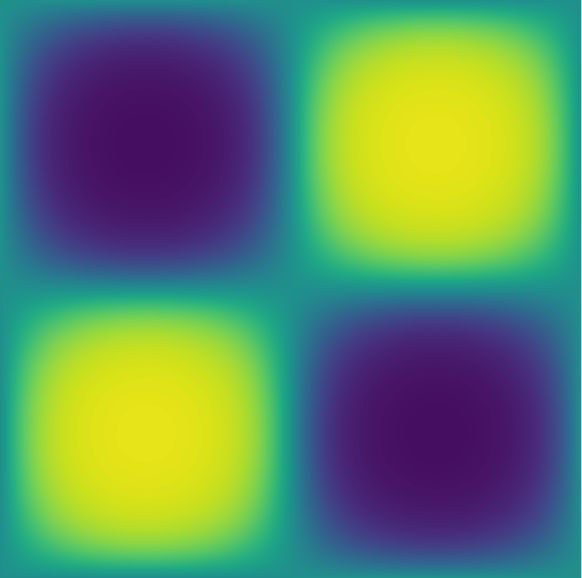}
        \hspace{0.7cm}
        \includegraphics[width=0.2\linewidth]{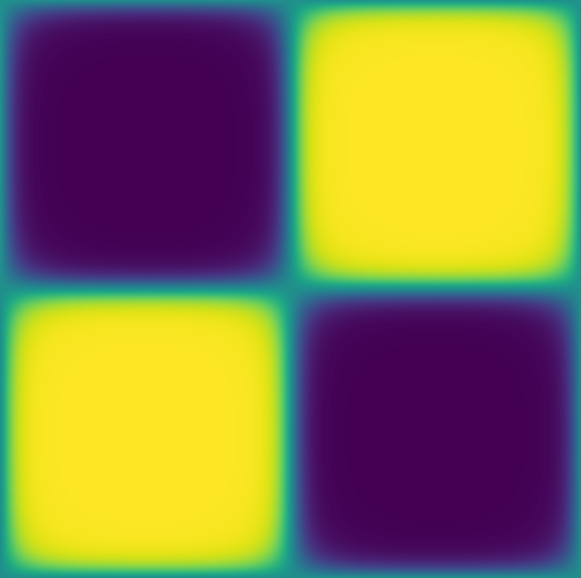}
        \hspace{0.7cm}
        \includegraphics[width=0.2\linewidth]{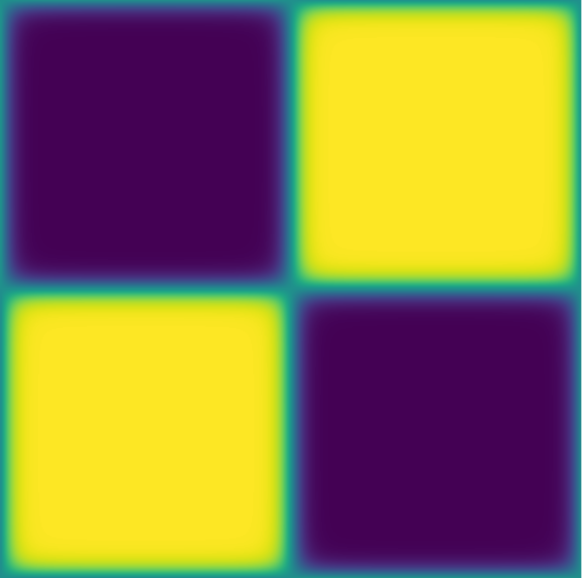}
        \caption{Allen--Cahn simulation, RKG2}
        \end{subfigure}
    \end{minipage}%
    \hspace{0.1cm}
    \begin{minipage}{0.06\linewidth}
        \vspace{-0.5cm}  
        \begin{tikzpicture}
            \begin{axis}[
                hide axis,
                scale only axis,
                height=6cm,
                width=0.1cm,
                colormap/viridis,
                colorbar right,
                colorbar style={
                    width=0.3cm,
                    ytick={-1, -0.5, 0, 0.5, 1},
                    ylabel={},
                },
                point meta min=-1,
                point meta max=1,
            ]
                \addplot [draw=none] coordinates {(0,0)};
            \end{axis}
        \end{tikzpicture}
    \end{minipage}
    \vspace{-0.5em}
    \caption{Allen--Cahn equation coarsening evolution at $t=2, 4, 6$ and $8$, with initial condition~\Cref{eqn:sinusoidal_ic_ac}.}
    \label{fig:ac_sinusoidal}
\end{figure}
\vspace{-1em}
\begin{figure}[!t]
    \centering
    \begin{subfigure}{0.4\linewidth}
    \captionsetup{justification=centering}
    \begin{overpic}[width=1\textwidth]{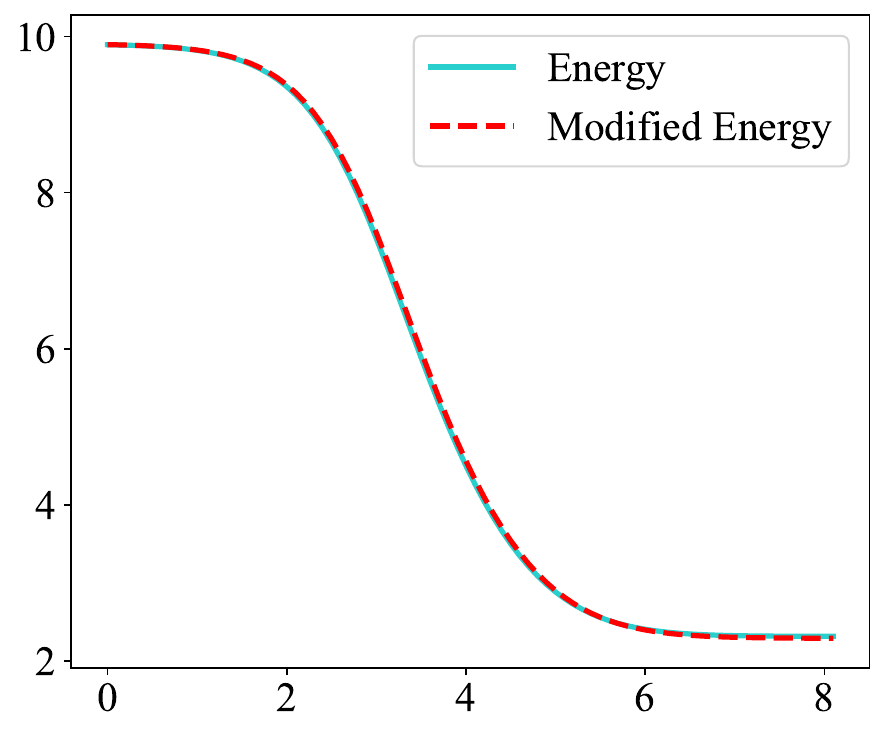}
        \put(52,-2){\small $t$}
        \put(-4,38){\rotatebox{90}{\small Energy}}
    \end{overpic}
    \caption{Allen--Cahn equation energy, RKL2}
    \label{subfig:allen_cahn_sinusoidal_energy_rkl2}
    \end{subfigure}
    \hspace{0.5em}
    \begin{subfigure}{0.4\linewidth}
    \captionsetup{justification=centering}
    \begin{overpic}[width=1\textwidth]{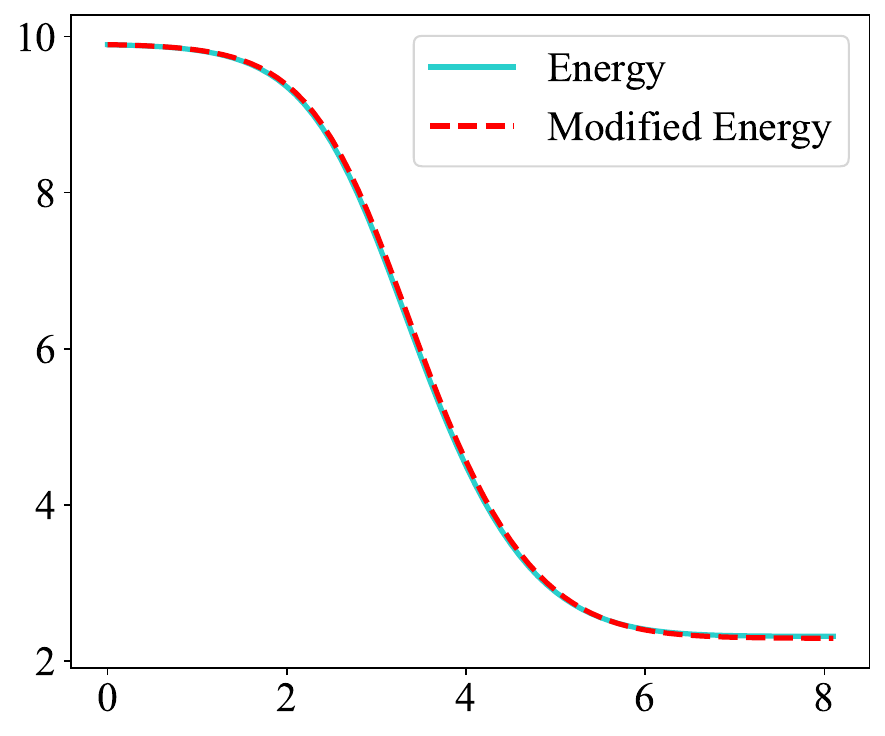}
        \put(52,-2){\small $t$}
        \put(-4,38){\rotatebox{90}{\small Energy}}
    \end{overpic}
    \caption{Allen--Cahn equation energy, RKG2}
    \label{subfig:allen_cahn_sinusoidal_energy_rkg2}
    \end{subfigure}
    \vspace{-0.5em}
    \captionsetup{width=0.9\linewidth}
    \caption{Original energy and modified energy dissipation of the Allen--Cahn equation simulations, with initial condition~\Cref{eqn:sinusoidal_ic_ac}. The modified energy is defined by~\Cref{eqn:nln_energy}.}
    \label{fig:ac_sinusoidal_energy}
\end{figure}

\subsection{Coarsening process and energy evolution in 2-dimensional space}
\label{exp:coarsening_2d_rectangular}
\begin{figure}[!t]
    \centering
    \begin{minipage}{0.8\linewidth}
        \begin{subfigure}{\linewidth}
        \captionsetup{justification=centering}
        \includegraphics[width=0.2\linewidth]{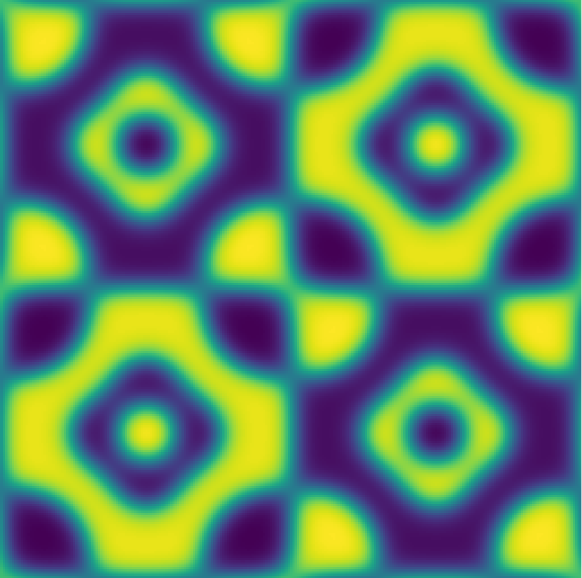}
        \hspace{0.7cm}
        \includegraphics[width=0.2\linewidth]{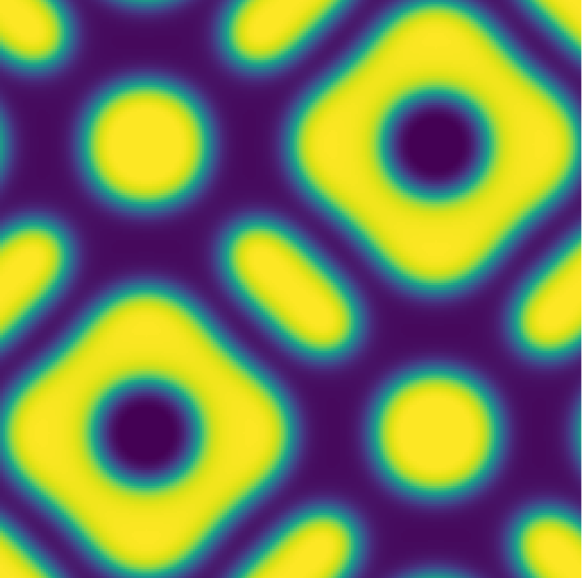}
        \hspace{0.7cm}
        \includegraphics[width=0.2\linewidth]{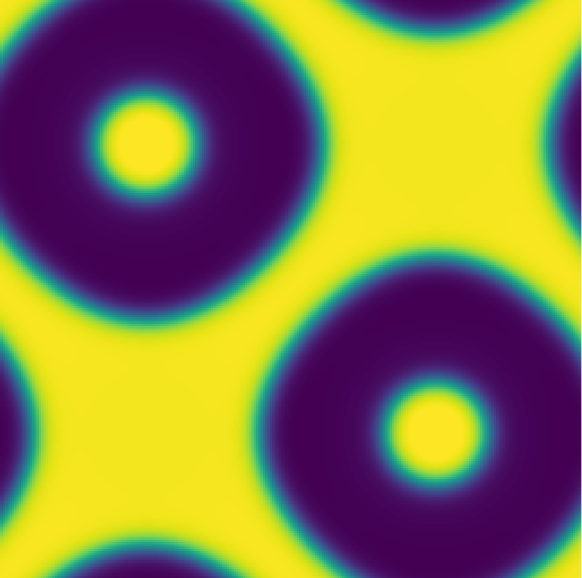}
        \hspace{0.7cm}
        \includegraphics[width=0.2\linewidth]{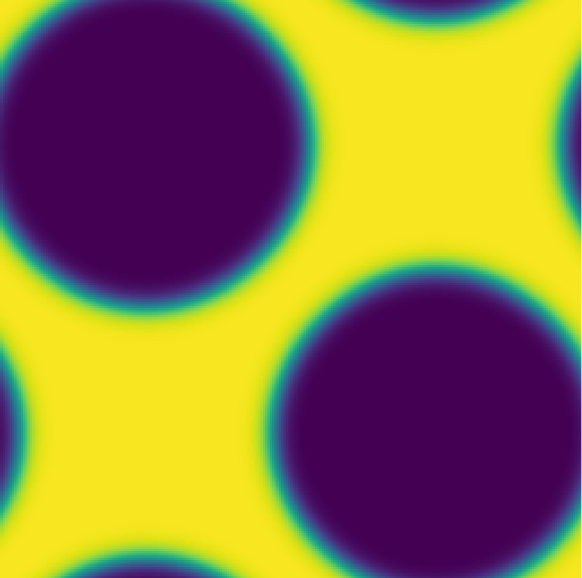}
        \caption{Cahn--Hilliard simulation, RKL2}
        \end{subfigure}
        \begin{subfigure}{\linewidth}
        \captionsetup{justification=centering}
        \includegraphics[width=0.2\linewidth]{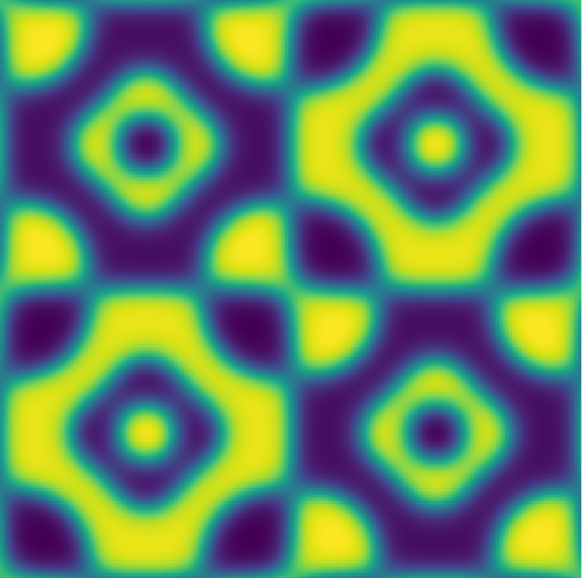}
        \hspace{0.7cm}
        \includegraphics[width=0.2\linewidth]{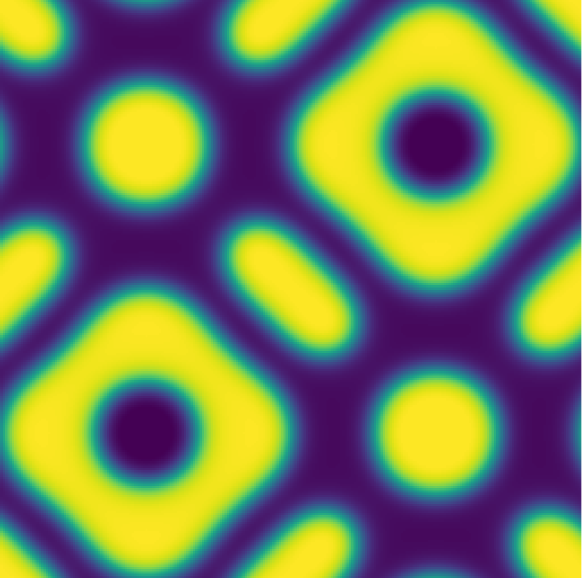}
        \hspace{0.7cm}
        \includegraphics[width=0.2\linewidth]{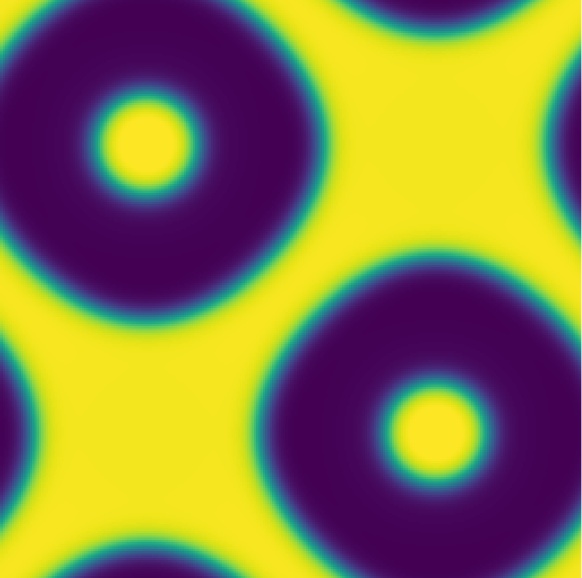}
        \hspace{0.7cm}
        \includegraphics[width=0.2\linewidth]{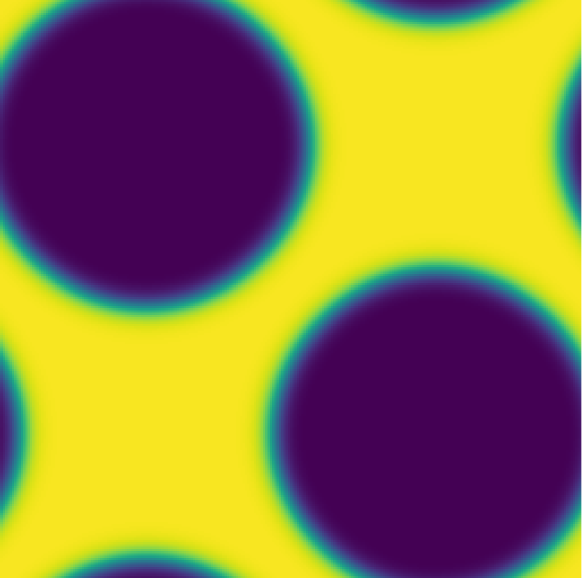}
        \caption{Cahn--Hilliard simulation, RKG2}
        \end{subfigure}
    \end{minipage}%
    \hspace{0.1cm}
    \begin{minipage}{0.06\linewidth}
        \vspace{-0.5cm}  
        \begin{tikzpicture}
            \begin{axis}[
                hide axis,
                scale only axis,
                height=6cm,
                width=0.1cm,
                colormap/viridis,   
                colorbar right,
                colorbar style={
                    width=0.3cm,
                    ytick={-1, -0.5, 0, 0.5, 1},   
                    ylabel={},
                },
                point meta min=-1,  
                point meta max=1,   
            ]
                \addplot [draw=none] coordinates {(0,0)};
            \end{axis}
        \end{tikzpicture}
    \end{minipage}
    \vspace{-0.5em}
    \caption{Cahn--Hilliard equation coarsening evolution at $t=1, 2, 5$ and $8$, with initial condition~\Cref{eqn:sinusoidal_ic_ch}.}
    \label{fig:ch_sinusoidal}
\end{figure}
\begin{figure}[!t]
    \centering
    \begin{subfigure}{0.4\linewidth}
    \captionsetup{justification=centering}
    \begin{overpic}[width=1\textwidth]{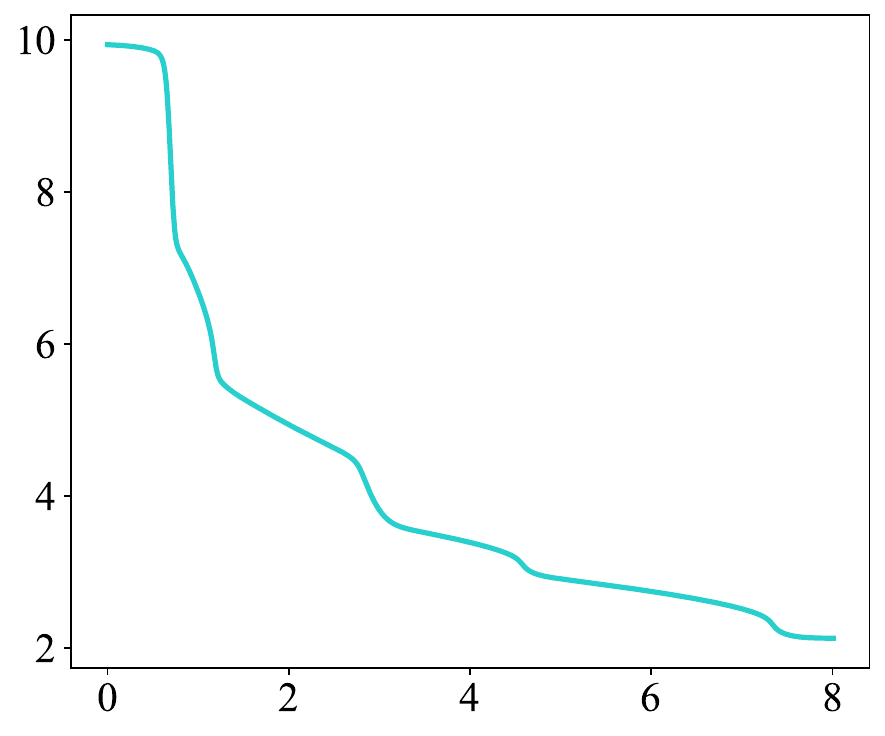}
        \put(52,-2){\small $t$}
        \put(-4,38){\rotatebox{90}{\small Energy}}
    \end{overpic}
    \caption{Cahn--Hilliard equation energy, RKL2}
    \label{subfig:cahn_hilliard_sinusoidal_energy_rkl2}
    \end{subfigure}
    \hspace{0.5em}
    \begin{subfigure}{0.4\linewidth}
    \captionsetup{justification=centering}
    \begin{overpic}[width=1\textwidth]{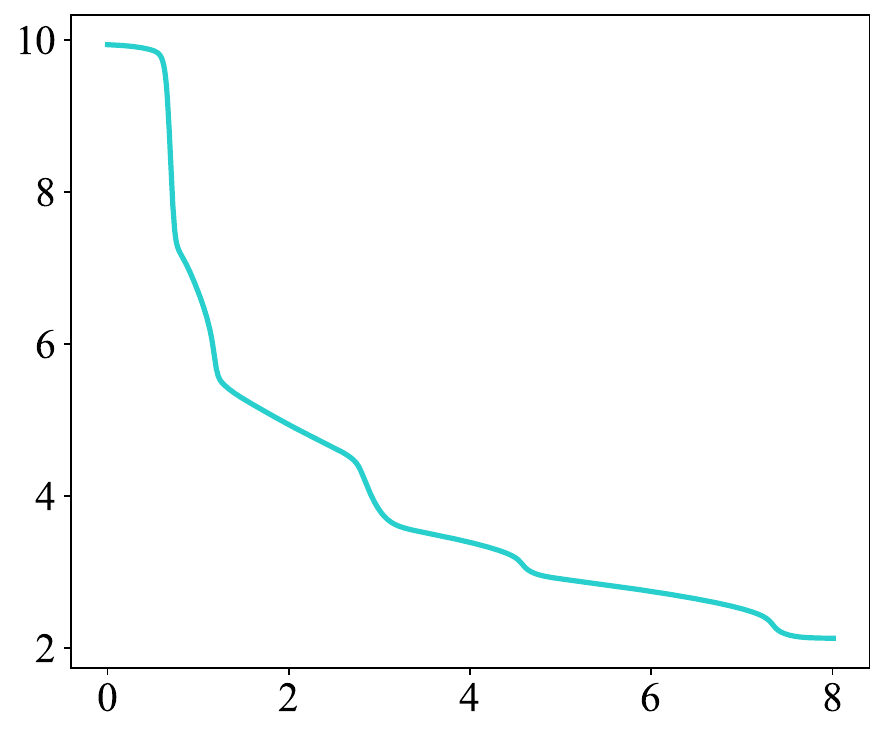}
        \put(52,-2){\small $t$}
        \put(-4,38){\rotatebox{90}{\small Energy}}
    \end{overpic}
    \caption{Cahn--Hilliard equation energy, RKG2}
    \label{subfig:cahn_hilliard_sinusoidal_energy_rkg2}
    \end{subfigure}
    \vspace{-0.5em}
    \caption{Energy dissipation of the Cahn--Hilliard equation simulations, with initial condition~\Cref{eqn:sinusoidal_ic_ch}.}
    \label{fig:ch_sinusoidal_energy}
\end{figure}
In this section, we simulate the rescaled Allen--Cahn and Cahn--Hilliard equation 
$$\frac{\partial u}{\partial t}=\epsilon^2\nabla^2u-F'(u), \frac{\partial u}{\partial t}=\nabla^2\left(-\epsilon^2\nabla^2u+F'(u)\right)$$
in the domain $[0, 2\pi]\times[0, 2\pi]$ with periodic boundary condition. We use $\epsilon=0.1$ for all the tests.
We first simulate the coarsening behavior to validate our method, by comparing the simulation results to the ones in the existing literature. For the Allen--Cahn equation, we use the initial condition
\begin{equation}
\label{eqn:sinusoidal_ic_ac}
u(x, y, 0)=0.05\sin(x)\sin(y).
\end{equation}
The spatial grid is 512$\times$512, with $\tau=1e$-$1$. We use $s=10$ for RKL2 and $s=12$ for the RKG2 methods so that their stability can be guaranteed, without an excessive number of substeps. We simulate until $T$=$8$. As shown in~\Cref{fig:ac_sinusoidal}, the results are visually indistinguishable between methods and are consistent with Figure 1 in~\cite{ETDRKEnergyStable}. The maximum principle is also preserved for the simulation. 

For the Cahn--Hilliard equation, it \delete{was}\revision{is} verified in~\cite{CHENG201545} that the symmetric sinusoidal initial condition above leads to an unstable evolution. Therefore, we instead use the asymmetric initial condition in~\cite{CHENG201545}
\begin{equation}
\label{eqn:sinusoidal_ic_ch}
u(x, y, 0)=0.05\sin(x)\sin(y)+0.001.
\end{equation}
The spatial grid is 256$\times$256. We use $s$=$200$ for RKL2 and $s$=$240$ for RKG2 to guarantee stability. The timestep is $\tau$=$1e$-$2$. We again simulate until $T$=$8$. The results in~\Cref{fig:ch_sinusoidal} agree with Figure 5.12 in~\cite{CHENG2020113070}, and using the fact that the energy and the solution do not change after about $T$=$8$. 

We also compute the energy for the numerical solutions; as shown in~\Cref{fig:ac_sinusoidal_energy} and~\Cref{fig:ch_sinusoidal_energy}, both the Allen--Cahn simulation and the Cahn--Hilliard simulations obtain monotonically decreasing energy curves; the shape and the values of the curves are also similar to the ones in Figure 2 of~\cite{ETDRK2MaximumPrinciple}. For the Allen--Cahn equation, we also compute the modified energy in~\Cref{eqn:nln_energy} in addition to the standard energy. For visualization purpose, when plotting the modified energy, we use the exact energy at $t=0$. As shown in~\Cref{fig:ac_sinusoidal_energy}, there is no visible difference between the exact energy and the modified energy.

\begin{figure}[!t]
    \centering
    \begin{minipage}{0.8\linewidth}
        \begin{subfigure}{\linewidth}
        \captionsetup{justification=centering}
        \includegraphics[width=0.2\linewidth]{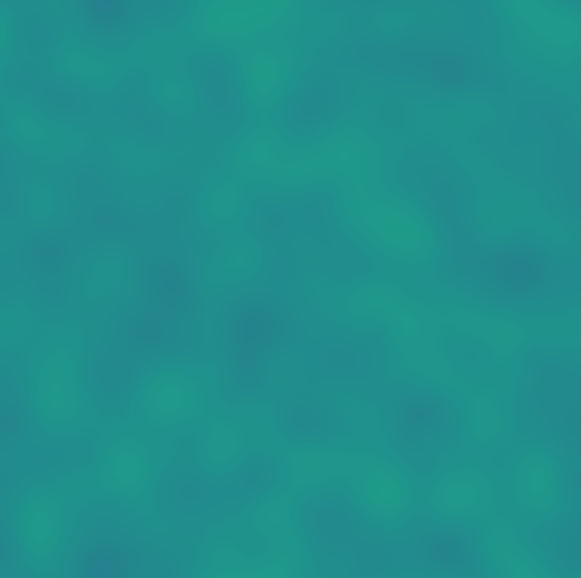}
        \hspace{0.7cm}
        \includegraphics[width=0.2\linewidth]{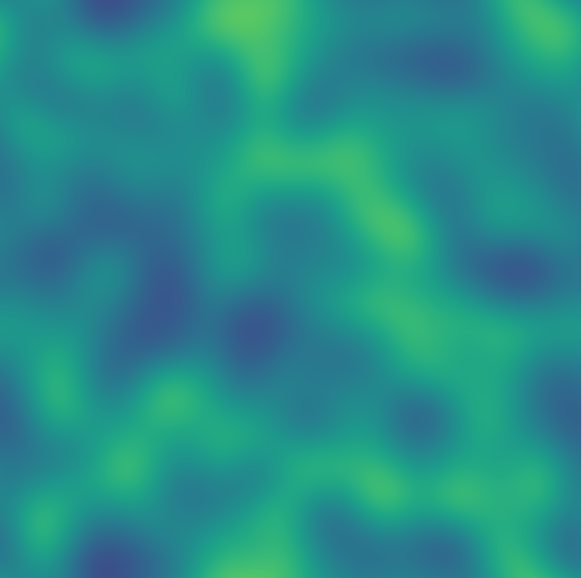}
        \hspace{0.7cm}
        \includegraphics[width=0.2\linewidth]{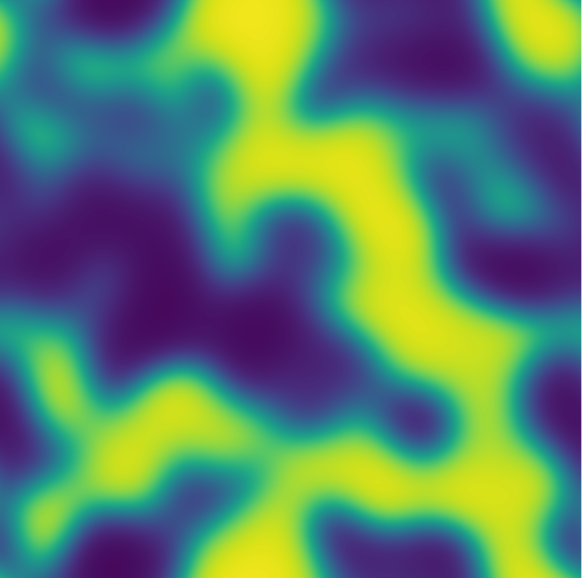}
        \hspace{0.7cm}
        \includegraphics[width=0.2\linewidth]{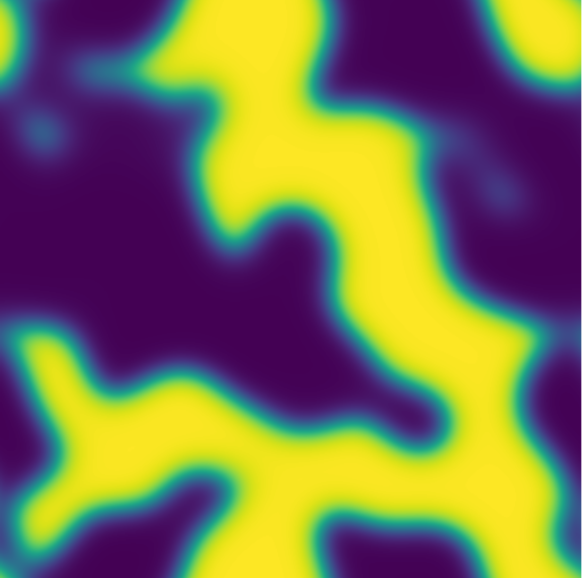}
        \caption{Allen--Cahn simulation, RKL2}
        \end{subfigure}
        \begin{subfigure}{\linewidth}
        \captionsetup{justification=centering}
        \includegraphics[width=0.2\linewidth]{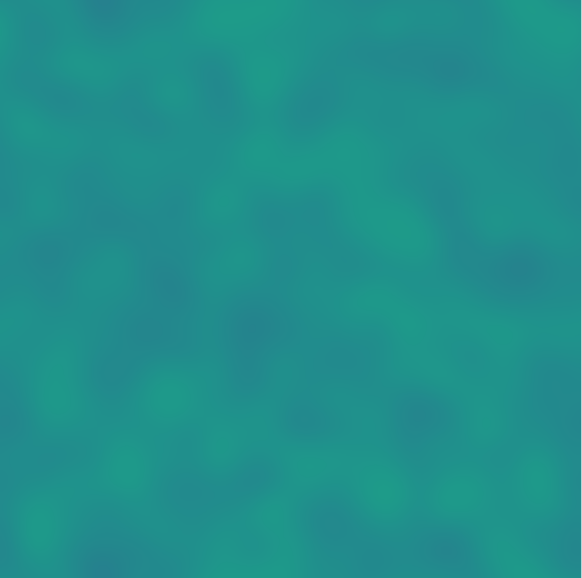}
        \hspace{0.7cm}
        \includegraphics[width=0.2\linewidth]{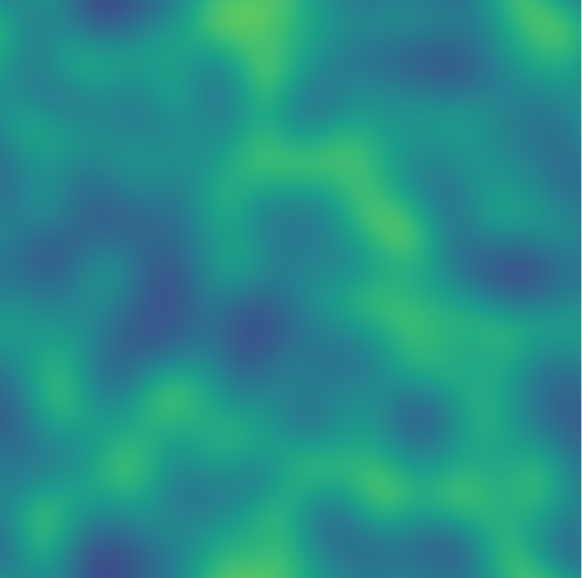}
        \hspace{0.7cm}
        \includegraphics[width=0.2\linewidth]{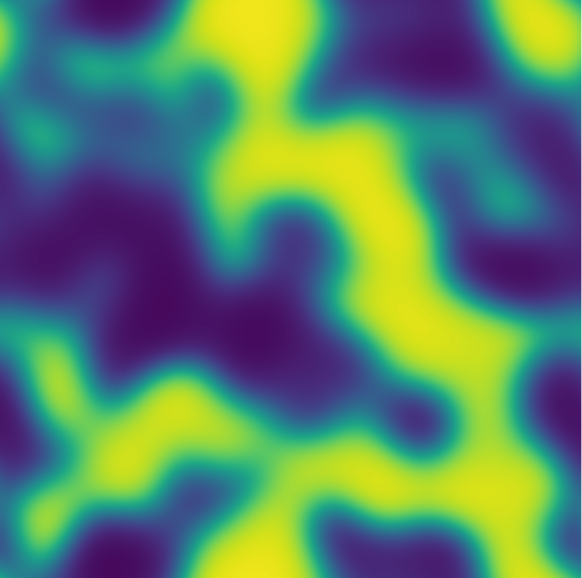}
        \hspace{0.7cm}
        \includegraphics[width=0.2\linewidth]{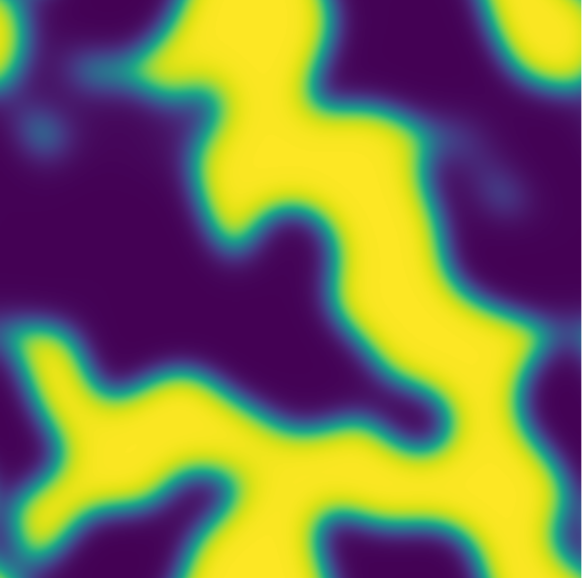}
        \caption{Allen--Cahn simulation, RKG2}
        \end{subfigure}
    \end{minipage}%
    \hspace{0.1cm}
    \begin{minipage}{0.06\linewidth}
        \vspace{-0.5cm}  
        \begin{tikzpicture}
            \begin{axis}[
                hide axis,
                scale only axis,
                height=6cm,
                width=0.1cm,
                colormap/viridis,   
                colorbar right,
                colorbar style={
                    width=0.3cm,
                    ytick={-1, -0.5, 0, 0.5, 1},   
                    ylabel={},
                },
                point meta min=-1,  
                point meta max=1,   
            ]
                \addplot [draw=none] coordinates {(0,0)};
            \end{axis}
        \end{tikzpicture}
    \end{minipage}
    \vspace{-0.5em}
    \caption{Allen--Cahn equation coarsening evolution at $t=2, 4, 6$ and $8$, with initial condition~\Cref{eqn:rand_ic_ac}.}
    \label{fig:ac_random}
\end{figure}
\begin{figure}[!t]
    \centering
    \begin{subfigure}{0.4\linewidth}
    \captionsetup{justification=centering}
    \begin{overpic}[width=1\textwidth]{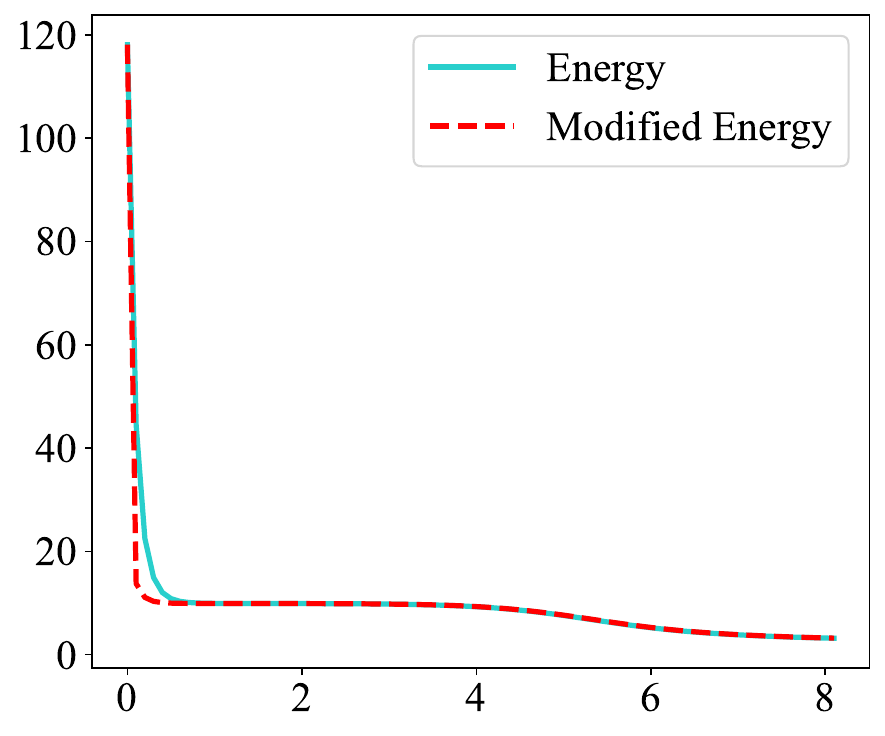}
        \put(52,-2){\small $t$}
        \put(-4,38){\rotatebox{90}{\small Energy}}
    \end{overpic}
    \caption{Allen--Cahn equation energy, RKL2}
    \end{subfigure}
    \hspace{0.5em}
    \begin{subfigure}{0.4\linewidth}
    \captionsetup{justification=centering}
    \begin{overpic}[width=1\textwidth]{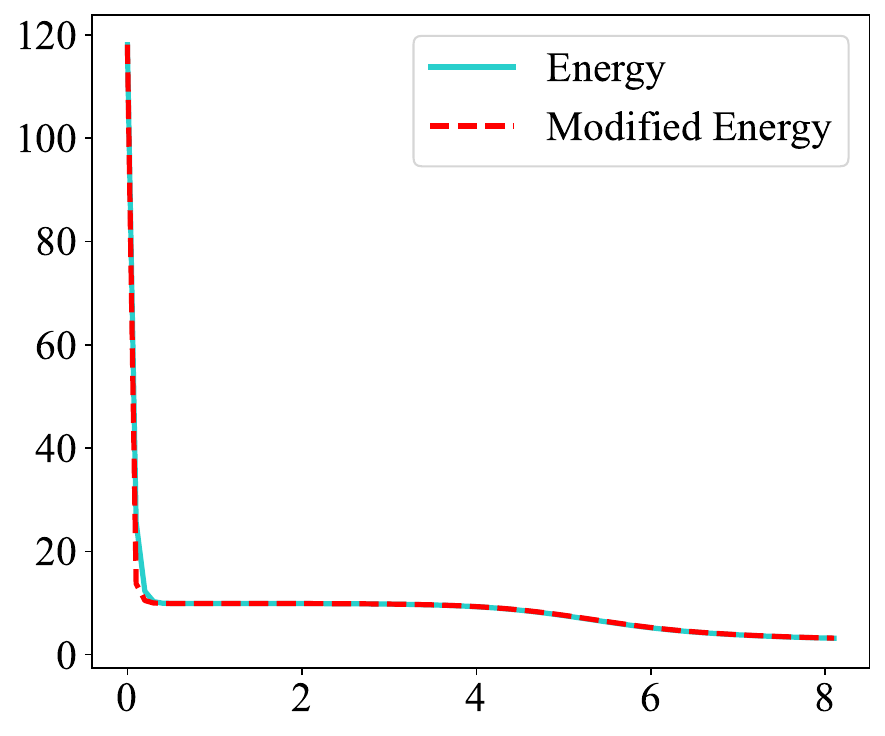}
        \put(52,-2){\small $t$}
        \put(-4,38){\rotatebox{90}{\small Energy}}
    \end{overpic}
    \caption{Allen--Cahn equation energy, RKG2}
    \end{subfigure}
    \vspace{-0.5em}
    \captionsetup{width=0.9\linewidth}
    \caption{Original energy and modified energy dissipation of the Allen--Cahn equation solutions, with initial condition~\Cref{eqn:rand_ic_ac}. The modified energy is defined by~\Cref{eqn:nln_energy}.}
    \label{fig:ac_random_energy}
\end{figure}

\begin{figure}[!t]
    \centering
    \begin{minipage}{0.8\linewidth}
        \begin{subfigure}{\linewidth}
        \captionsetup{justification=centering}
        \includegraphics[width=0.2\linewidth]{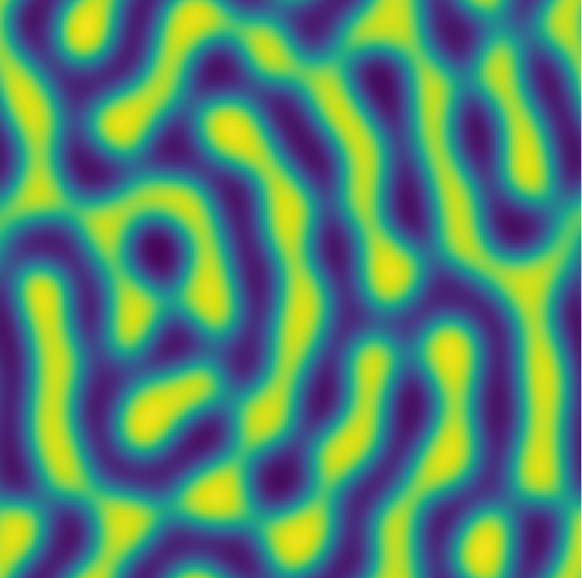}
        \hspace{0.7cm}
        \includegraphics[width=0.2\linewidth]{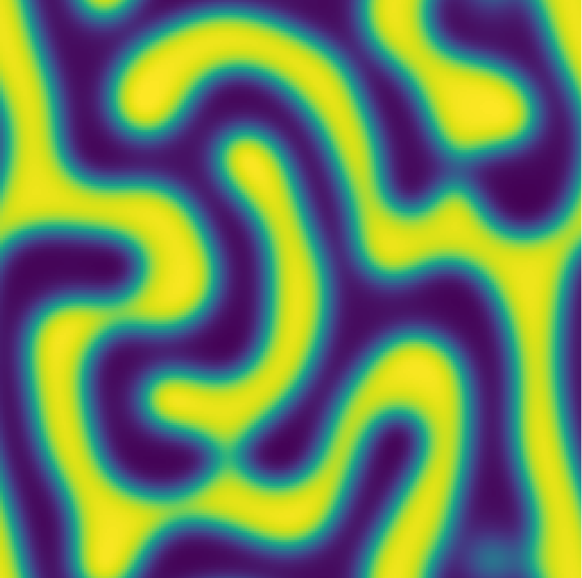}
        \hspace{0.7cm}
        \includegraphics[width=0.2\linewidth]{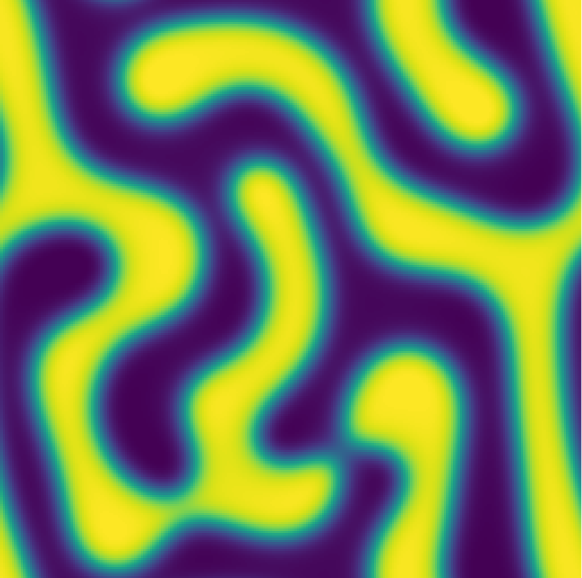}
        \hspace{0.7cm}
        \includegraphics[width=0.2\linewidth]{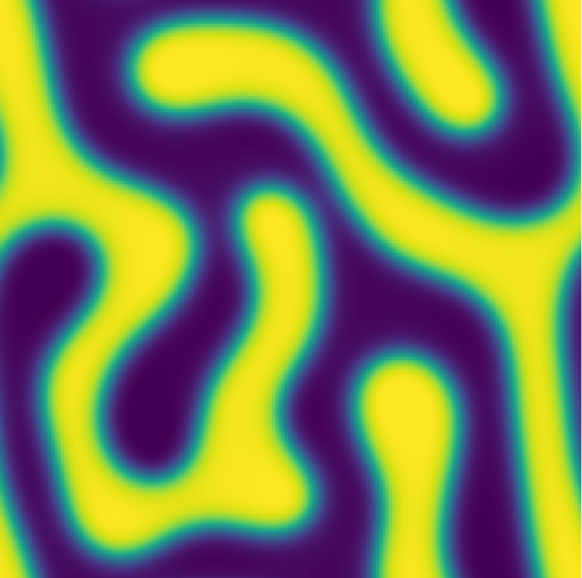}
        \caption{Cahn--Hilliard simulation, RKL2}
        \end{subfigure}
        \begin{subfigure}{\linewidth}
        \captionsetup{justification=centering}
        \includegraphics[width=0.2\linewidth]{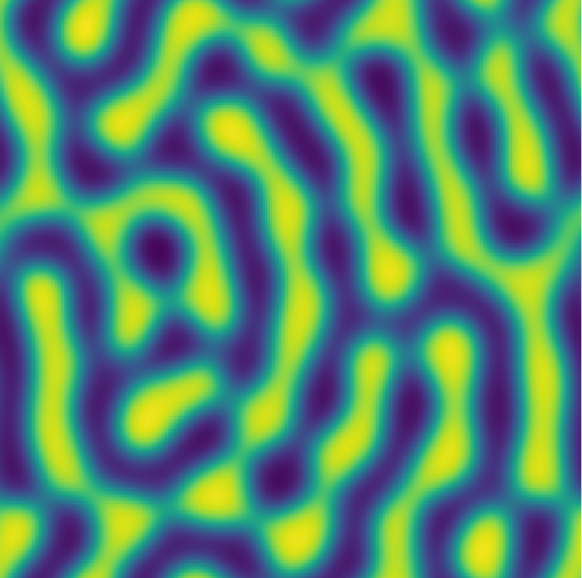}
        \hspace{0.7cm}
        \includegraphics[width=0.2\linewidth]{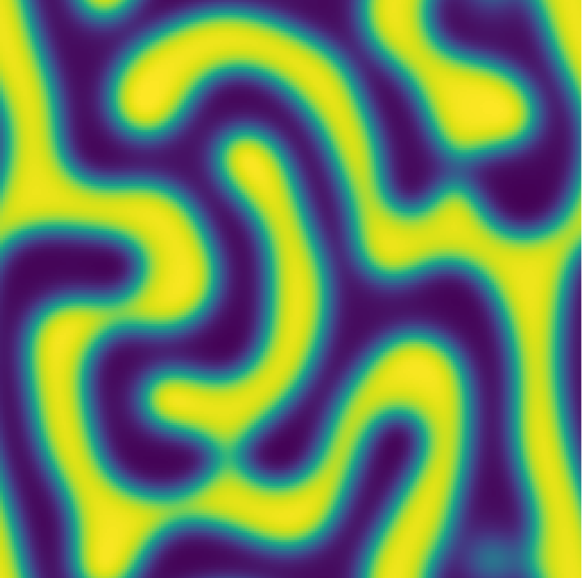}
        \hspace{0.7cm}
        \includegraphics[width=0.2\linewidth]{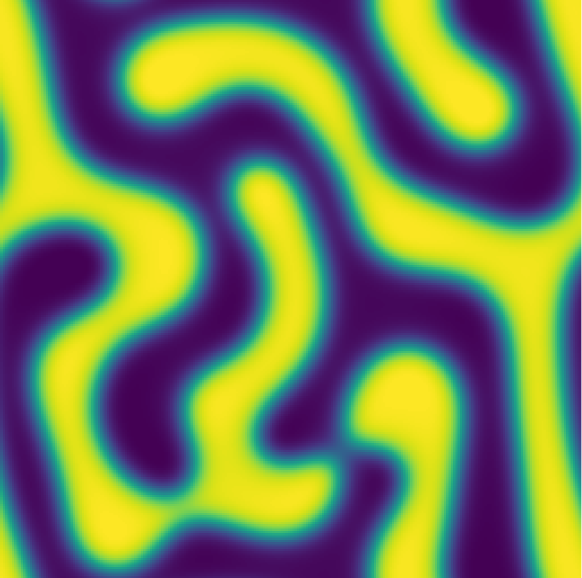}
        \hspace{0.7cm}
        \includegraphics[width=0.2\linewidth]{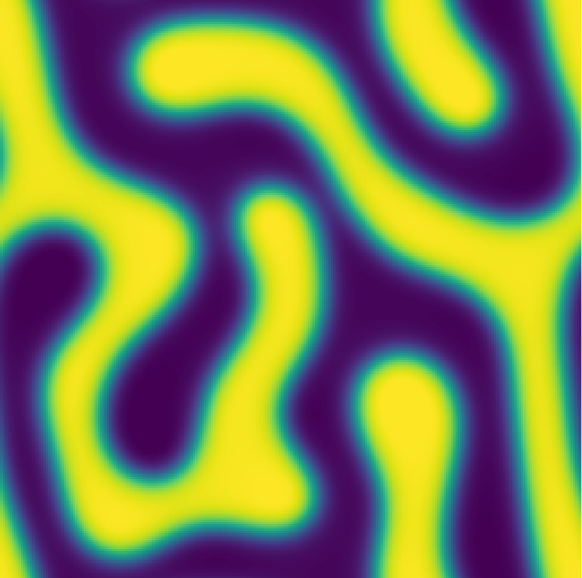}
        \caption{Cahn--Hilliard simulation, RKG2}
        \end{subfigure}
    \end{minipage}%
    \hspace{0.1cm}%
    \begin{minipage}{0.06\linewidth}
        \vspace{-0.5cm}
        \begin{tikzpicture}
            \begin{axis}[
                hide axis,
                scale only axis,
                height=6cm,
                width=0.1cm,
                colormap/viridis,
                colorbar right,
                colorbar style={
                    width=0.3cm,
                    ytick={-1, -0.5, 0, 0.5, 1},
                    ylabel={},
                },
                point meta min=-1,
                point meta max=1,
            ]
                \addplot [draw=none] coordinates {(0,0)};
            \end{axis}
        \end{tikzpicture}
    \end{minipage}
    \vspace{-1em}
    \caption{Cahn--Hilliard equation coarsening evolution at $t=0.5, 1, 1.5$ and $2$, with initial condition~\Cref{eqn:rand_ic_ch}.}
    \label{fig:ch_random}
\end{figure}
\begin{figure}[!t]
    \centering
    \begin{subfigure}{0.4\linewidth}
    \captionsetup{justification=centering}
    \begin{overpic}[width=1\textwidth]{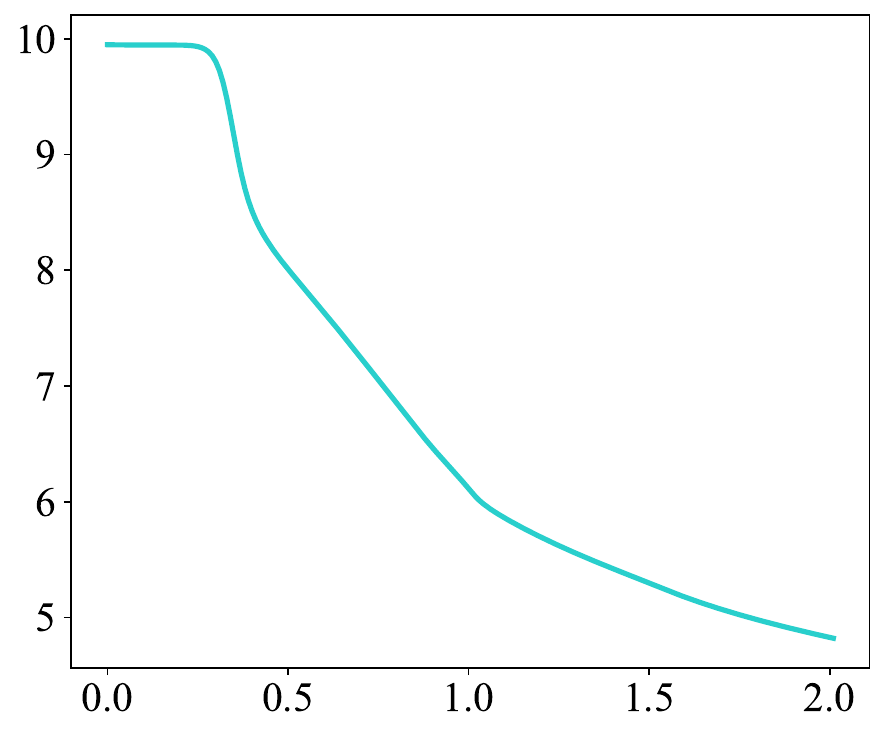}
        \put(52,-2){\small $t$}
        \put(-4,38){\rotatebox{90}{\small Energy}}
    \end{overpic}
    \caption{Cahn--Hilliard equation energy, RKL2}
    \end{subfigure}
    \hspace{0.5em}
    \begin{subfigure}{0.4\linewidth}
    \captionsetup{justification=centering}
    \begin{overpic}[width=1\textwidth]{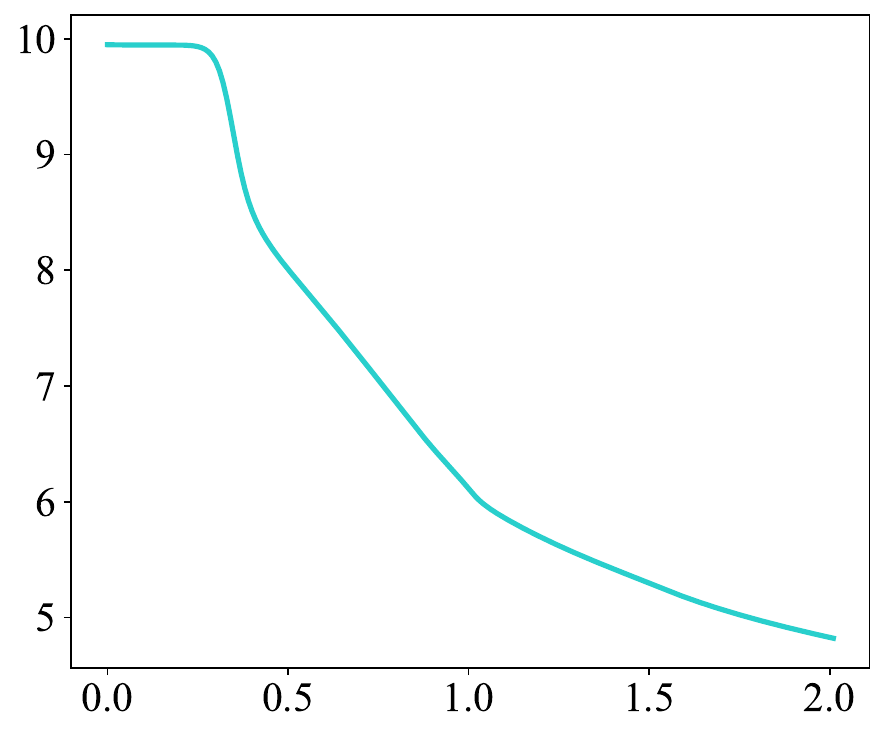}
        \put(52,-2){\small $t$}
        \put(-4,38){\rotatebox{90}{\small Energy}}
    \end{overpic}
    \caption{Cahn--Hilliard equation energy, RKG2}
    \end{subfigure}
    \caption{Energy dissipation of the Cahn--Hilliard equation simulations, with initial condition~\Cref{eqn:rand_ic_ch}.}
    \label{fig:ch_random_energy}
\end{figure}

Then we simulate Allen--Cahn and Cahn--Hilliard equations with random initial conditions. For Allen--Cahn, we use the initial condition
\begin{equation}
\label{eqn:rand_ic_ac}
u(x, y, 0)=0.5\text{Rand}(x, y),
\end{equation}
with the same settings as above. For Cahn--Hilliard, we again use an asymmetric random initial condition
\begin{equation}
\label{eqn:rand_ic_ch}
u(x, y, 0)=0.0025+0.004\text{Rand}(x, y),
\end{equation}
where $\text{Rand}(x, y)$ denotes random numbers sampled from $[-1, 1]$ uniformly. For Allen--Cahn, the spatial discretization is chosen to be 512$\times$512, the timestep is chosen to be $\tau$=$1e$-$1$, and we simulate until $T=8$. We still use $s$=$10$ for the RKL2 method and $s$=$12$ for the RKG2 method. For Cahn--Hilliard, the spatial discretization is chosen to be 256$\times$256, the timstep is chosen to be $\tau$=$1e$-$2$, and we simulate until $T=2$. The coarsening parameter is chosen to be $\epsilon=0.1$. We again use $s$=$200$ for RKL2 and $s$=$240$ for RKG2 as in the previous example. Periodic boundary condition is applied. The coarsening behaviors are shown in~\Cref{fig:ac_random} and~\Cref{fig:ch_random}; as shown, both RKL2 and RKG2 produce visually the same coarsening behaviors. \revision{Notice that for consistency, we use the color bar with range $[-1, 1]$ for the Cahn--Hilliard equation plots, although some of the simulation values might exceed this range with trivial amount due to the lack of maximum principle; visually the effect is negligible.} The energy curves of the methods are shown in~\Cref{fig:ac_random_energy} and~\Cref{fig:ch_random_energy}, both of which are monotonically decreasing. The Allen--Cahn simulations in addition have almost identical energy and modified energy as shown in~\Cref{fig:ac_random_energy}.

\subsection{Coarsening process and energy evolution in 2-dimensional irregular domain}
\begin{figure}[!t]
    \centering
    \begin{minipage}{0.8\linewidth}
        \begin{subfigure}{\linewidth}
        \captionsetup{justification=centering}
        \includegraphics[width=0.2\linewidth]{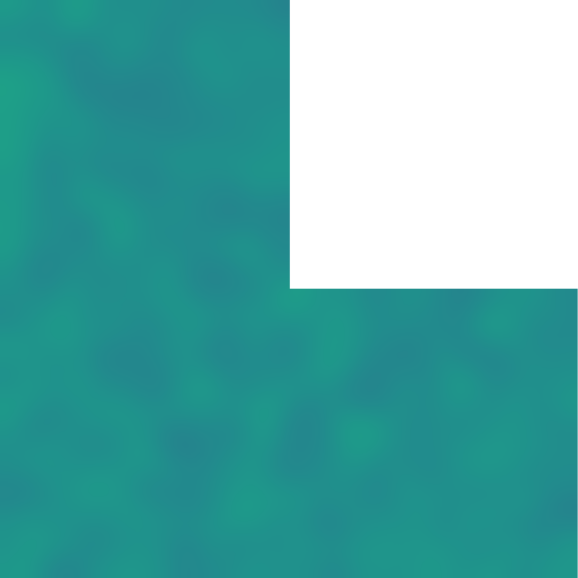}
        \hspace{0.7cm}
        \includegraphics[width=0.2\linewidth]{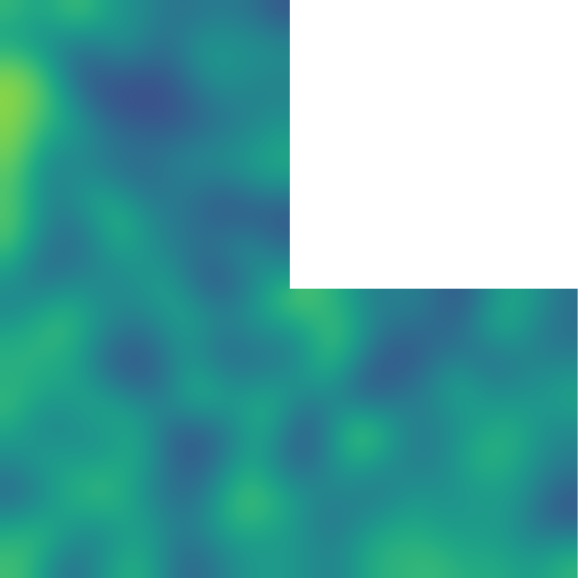}
        \hspace{0.7cm}
        \includegraphics[width=0.2\linewidth]{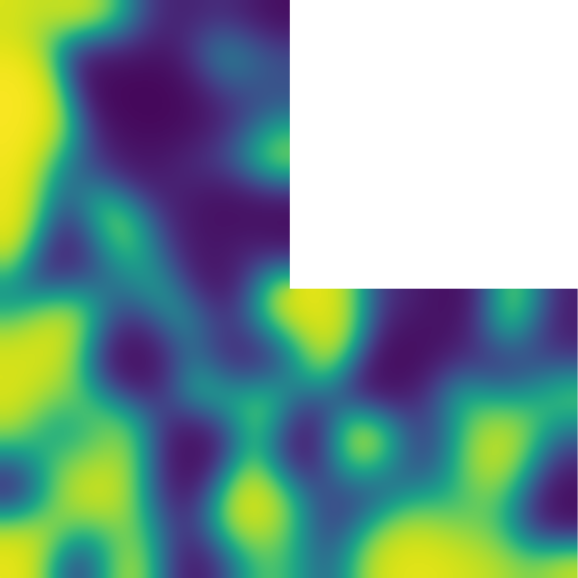}
        \hspace{0.7cm}
        \includegraphics[width=0.2\linewidth]{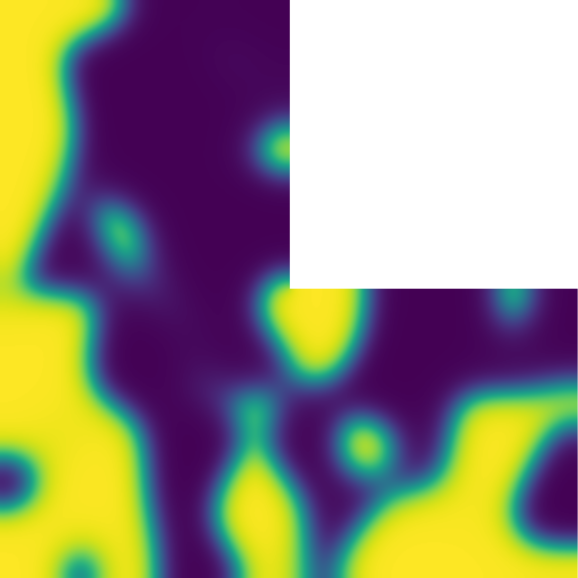}
        \caption{Allen--Cahn simulation, RKL2}
        \end{subfigure}
        \begin{subfigure}{\linewidth}
        \captionsetup{justification=centering}
        \includegraphics[width=0.2\linewidth]{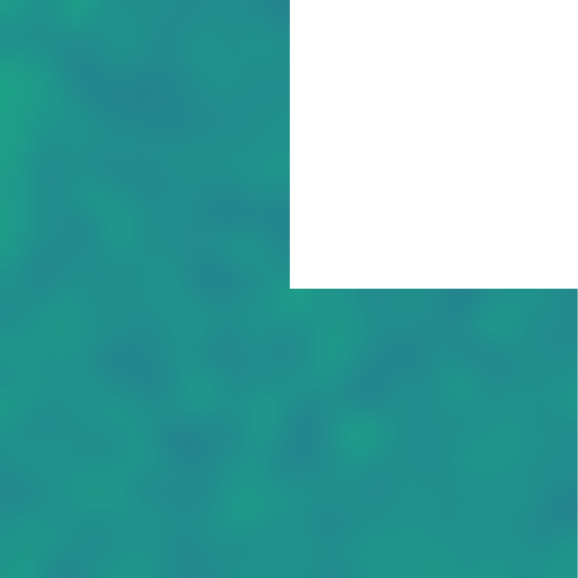}
        \hspace{0.7cm}
        \includegraphics[width=0.2\linewidth]{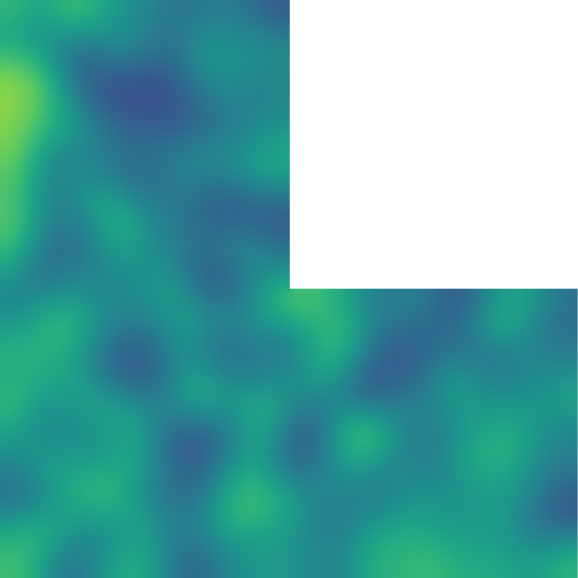}
        \hspace{0.7cm}
        \includegraphics[width=0.2\linewidth]{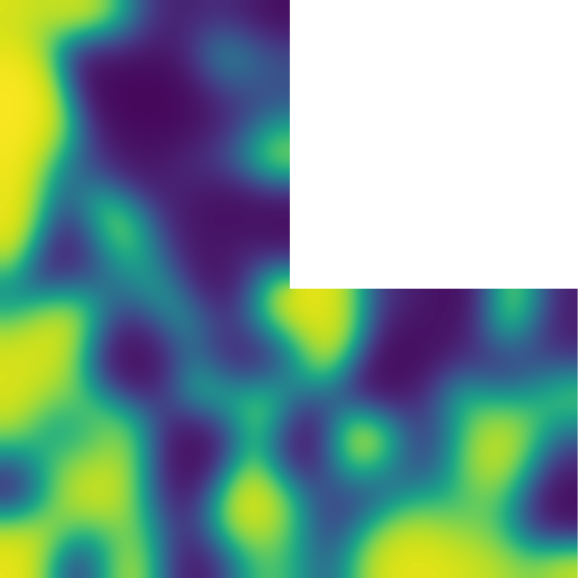}
        \hspace{0.7cm}
        \includegraphics[width=0.2\linewidth]{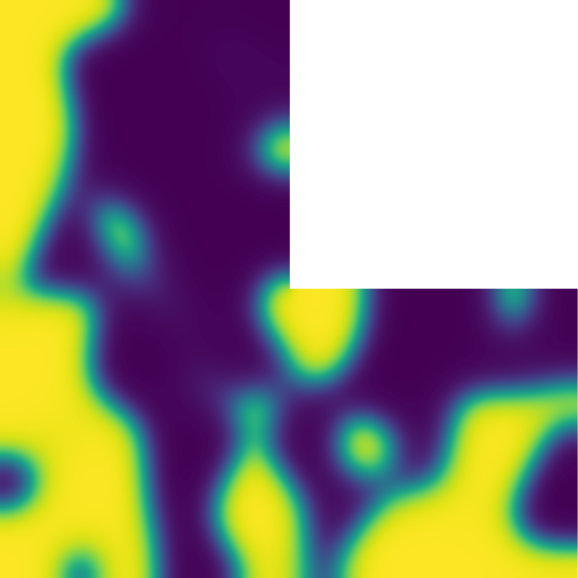}
        \caption{Allen--Cahn simulation, RKG2}
        \end{subfigure}
    \end{minipage}%
    \hspace{0.1cm}%
    \begin{minipage}{0.06\linewidth}
        \vspace{-0.5cm}
        \begin{tikzpicture}
            \begin{axis}[
                hide axis,
                scale only axis,
                height=6cm,
                width=0.1cm,
                colormap/viridis,
                colorbar right,
                colorbar style={
                    width=0.3cm,
                    ytick={-1, -0.5, 0, 0.5, 1},
                    ylabel={},
                },
                point meta min=-1,
                point meta max=1,
            ]
                \addplot [draw=none] coordinates {(0,0)};
            \end{axis}
        \end{tikzpicture}
    \end{minipage}
    \vspace{-1em}
    \caption{Allen--Cahn equation coarsening evolution at $t=2, 4, 6$ and $8$, within the $L$-shaped irregular domain and initial condition~\Cref{eqn:rand_ic_ac}.}
    \label{fig:ac_random_irregular}
\end{figure}
\begin{figure}[!t]
    \centering
    \begin{subfigure}{0.4\linewidth}
    \captionsetup{justification=centering}
    \begin{overpic}[width=1\textwidth]{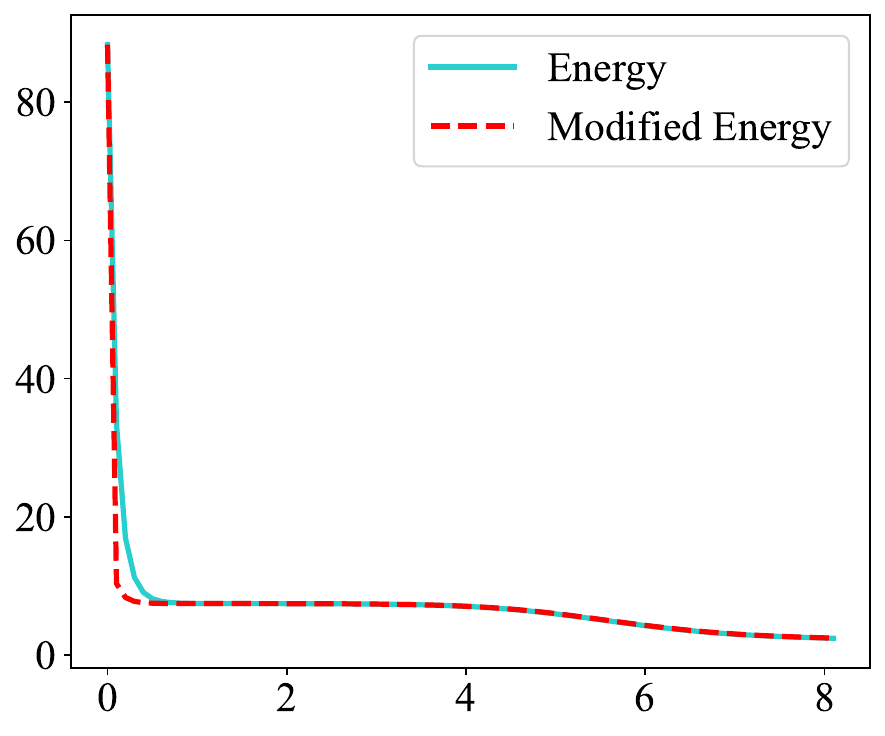}
        \put(52,-2){\small $t$}
        \put(-4,38){\rotatebox{90}{\small Energy}}
    \end{overpic}
    \caption{Allen--Cahn equation energy, RKL2}
    \end{subfigure}
    \hspace{0.5em}
    \begin{subfigure}{0.4\linewidth}
    \captionsetup{justification=centering}
    \begin{overpic}[width=1\textwidth]{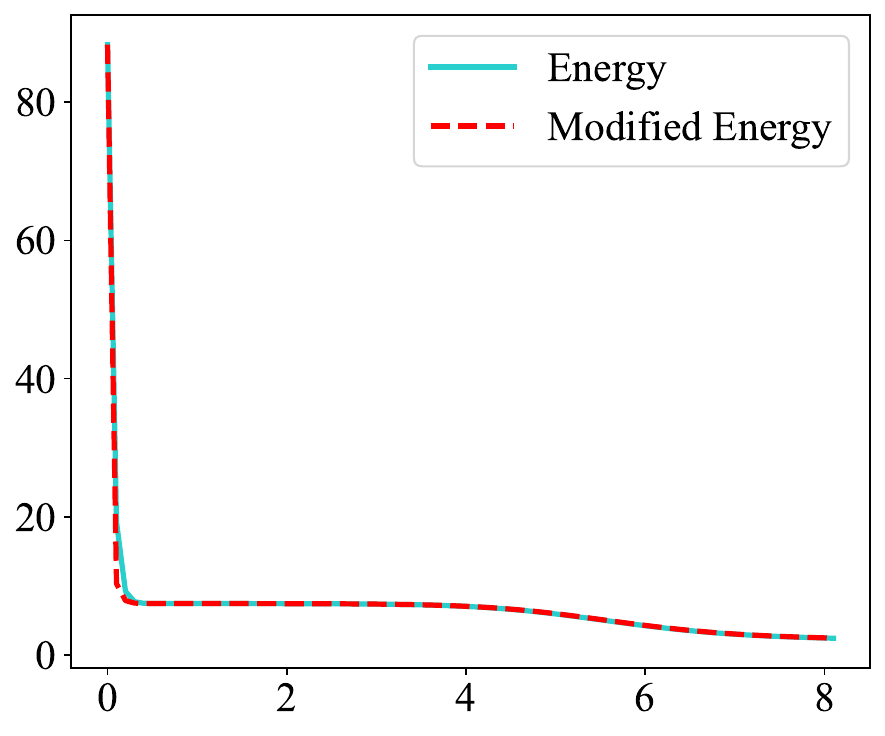}
        \put(52,-2){\small $t$}
        \put(-4,38){\rotatebox{90}{\small Energy}}
    \end{overpic}
    \caption{Allen--Cahn equation energy, RKG2}
    \end{subfigure}
    \vspace{-0.5em}
    \captionsetup{width=0.9\linewidth}
    \caption{Energy dissipation of the Allen--Cahn equation simulations on $L$-shaped domain, with initial condition~\Cref{eqn:rand_ic_ac}. The modified energy is defined by~\Cref{eqn:nln_energy}.}
    \label{fig:ac_random_energy_irregular}
\end{figure}

\begin{figure}[!t]
    \centering
    \begin{minipage}{0.8\linewidth}
        \begin{subfigure}{\linewidth}
        \captionsetup{justification=centering}
        \includegraphics[width=0.2\linewidth]{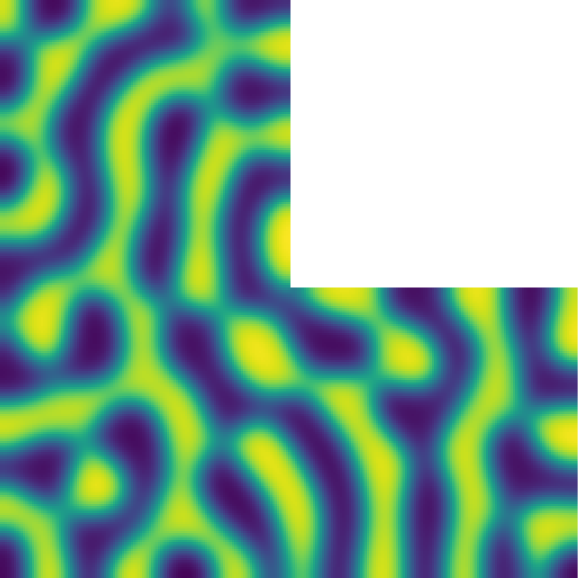}
        \hspace{0.7cm}
        \includegraphics[width=0.2\linewidth]{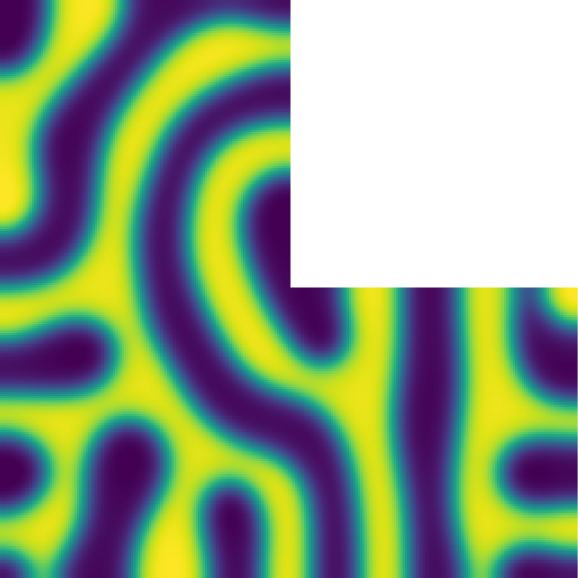}
        \hspace{0.7cm}
        \includegraphics[width=0.2\linewidth]{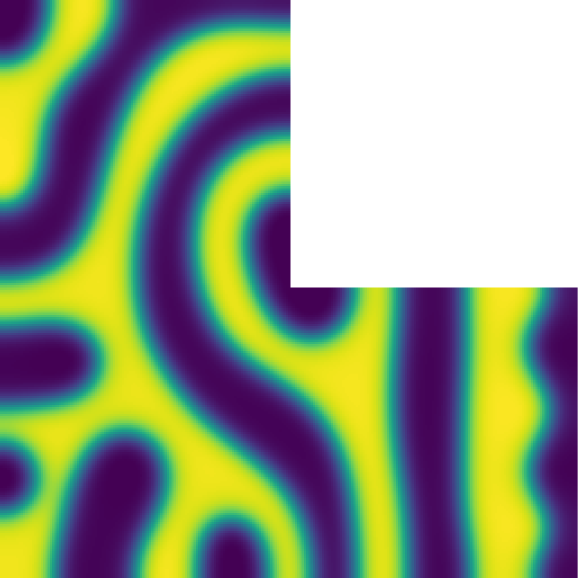}
        \hspace{0.7cm}
        \includegraphics[width=0.2\linewidth]{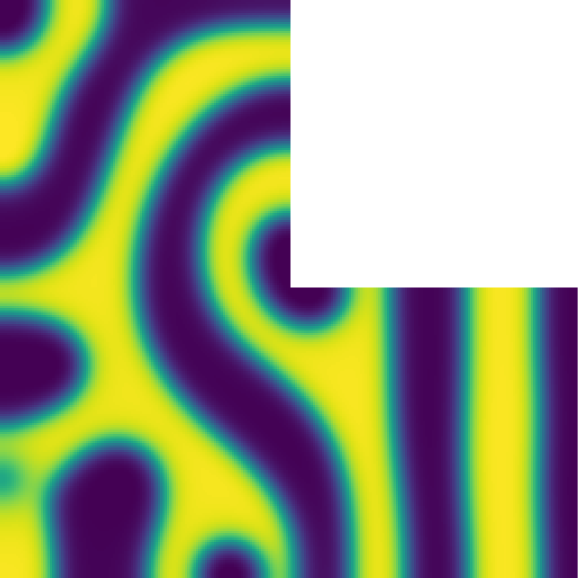}
        \caption{Cahn--Hilliard simulation on $L$-shaped domain, RKL2}
        \end{subfigure}
        \begin{subfigure}{\linewidth}
        \captionsetup{justification=centering}
        \includegraphics[width=0.2\linewidth]{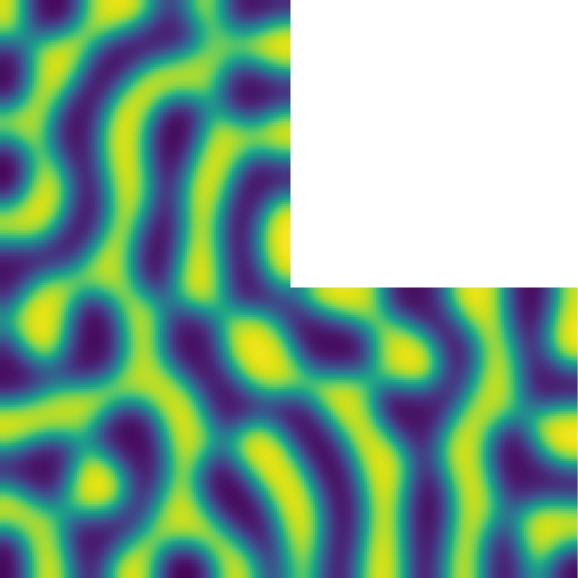}
        \hspace{0.7cm}
        \includegraphics[width=0.2\linewidth]{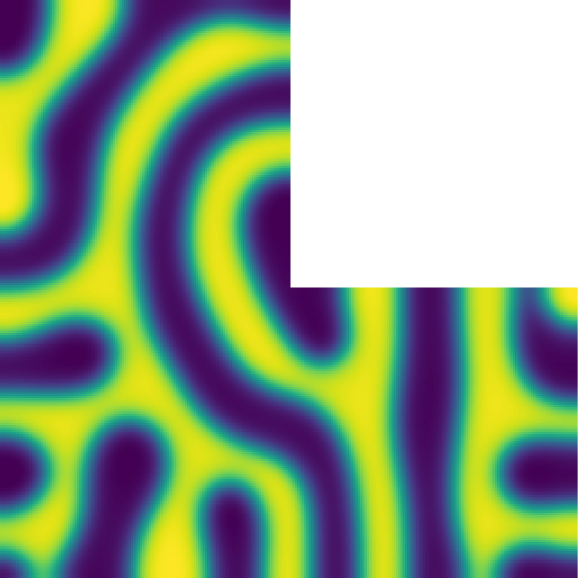}
        \hspace{0.7cm}
        \includegraphics[width=0.2\linewidth]{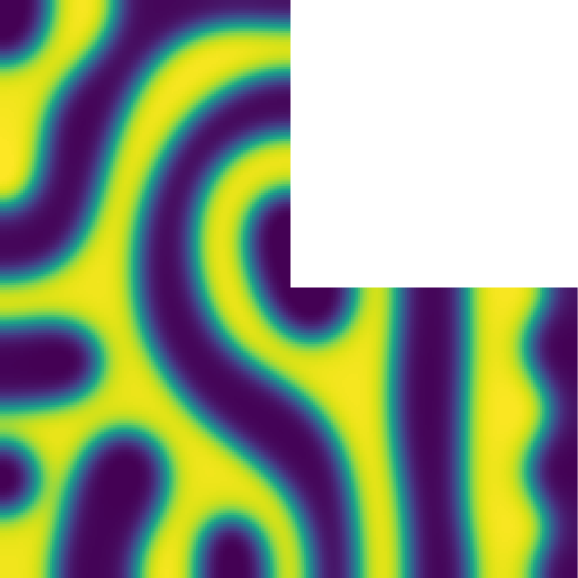}
        \hspace{0.7cm}
        \includegraphics[width=0.2\linewidth]{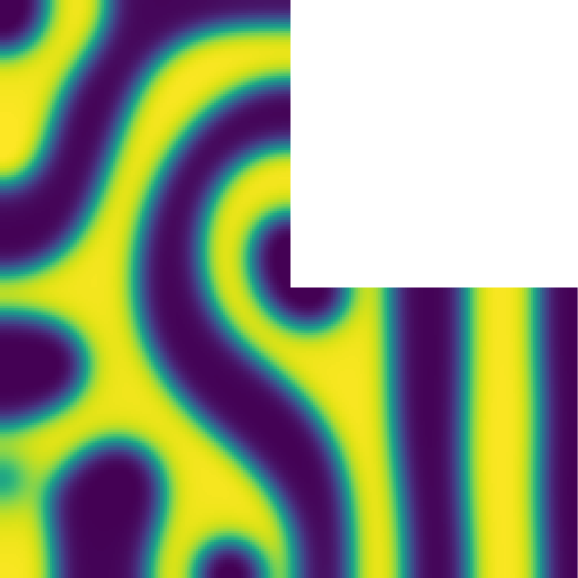}
        \caption{Cahn--Hilliard simulation on $L$-shaped domain, RKG2}
        \end{subfigure}
    \end{minipage}%
    \hspace{0.1cm}%
    \begin{minipage}{0.06\linewidth}
        \vspace{-0.5cm}
        \begin{tikzpicture}
            \begin{axis}[
                hide axis,
                scale only axis,
                height=6cm,
                width=0.1cm,
                colormap/viridis,
                colorbar right,
                colorbar style={
                    width=0.3cm,
                    ytick={-1, -0.5, 0, 0.5, 1},
                    ylabel={},
                },
                point meta min=-1,
                point meta max=1,
            ]
                \addplot [draw=none] coordinates {(0,0)};
            \end{axis}
        \end{tikzpicture}
    \end{minipage}
    \vspace{-1em}
    \caption{Cahn--Hilliard equation coarsening evolution at $t=0.5, 1, 1.5$ and $2$, within the $L$-shaped irregular domain and initial condition~\Cref{eqn:rand_ic_ch}. As shown, there is no visual difference between the RKL2 and RKG2 simulations.}
    \label{fig:ch_random_irregular}
\end{figure}

\begin{figure}[!t]
    \centering
    \begin{subfigure}{0.4\linewidth}
    \captionsetup{justification=centering}
    \begin{overpic}[width=1\textwidth]{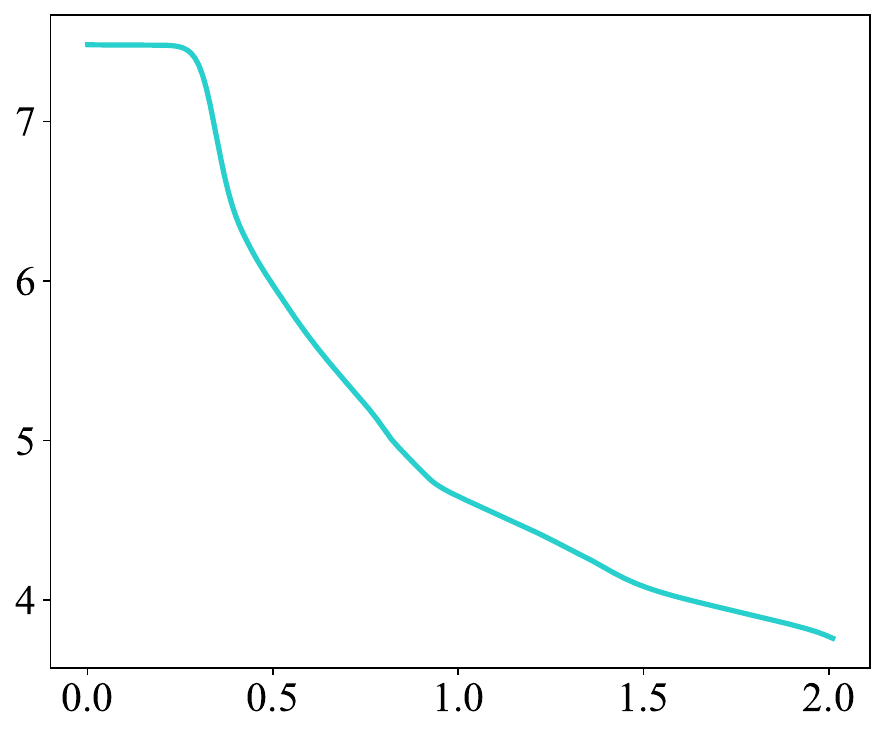}
        \put(52,-2){\small $t$}
        \put(-4,38){\rotatebox{90}{\small Energy}}
    \end{overpic}
    \caption{Cahn--Hilliard equation energy, RKL2}
    \end{subfigure}
    \hspace{0.5em}
    \begin{subfigure}{0.4\linewidth}
    \captionsetup{justification=centering}
    \begin{overpic}[width=1\textwidth]{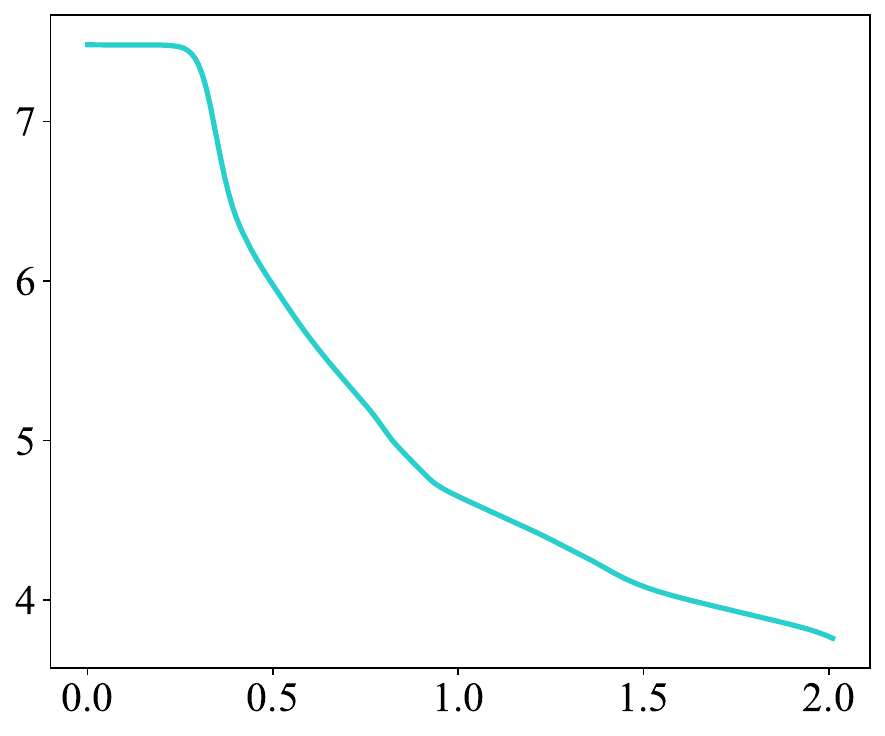}
        \put(52,-2){\small $t$}
        \put(-4,38){\rotatebox{90}{\small Energy}}
    \end{overpic}
    \caption{Cahn--Hilliard equation energy, RKG2}
    \end{subfigure}
    \caption{Energy dissipation of the Cahn--Hilliard equation simulations on $L$-shaped domain, with initial condition~\Cref{eqn:rand_ic_ch}.}
    \label{fig:ch_random_energy_irregular}
\end{figure}

As discussed in~\Cref{sec:literature_review}, the exponential time differencing Runge--Kutta (ETDRK) methods in rectangular domains with periodic boundary condition can be solved extremely efficiently, using the Fourier spectral method for spatial discretization. However, it is also known that for irregular domains with non-periodic boundary conditions, Fourier spectral method cannot be applied, and other methods such as Chebyshev or Legendre spectral methods need to be utilized in general~\cite{spectralmethodtextbook}. Therefore, RKL and RKG methods are attractive as explicit methods that do not rely on the Fourier spectral approximation. We hence conduct numerical simulations on an $L$-shaped domain for Allen--Cahn and Cahn--Hilliard equations, using RKL and RKG methods; the domain is chosen to be $[0, 2\pi]^2\setminus[\pi, 2\pi]^2$. We again choose the random initial condition~\Cref{eqn:rand_ic_ac} and~\Cref{eqn:rand_ic_ch}, with the same parameters as in the case of rectangular domain. The boundary condition is chosen to be homogeneous Neumann boundary condition. As shown in~\Cref{fig:ac_random_irregular},~\Cref{fig:ch_random_irregular},~\Cref{fig:ac_random_energy_irregular}, and~\Cref{fig:ch_random_energy_irregular}, the simulations of both methods are visually the same and the energy curves are decreasing. More importantly, the modified energy and the original energy again perfectly match each other, which confirms our claim that the proof of~\Cref{thm:nln_energy} does not rely on the fact that the domain is rectangular.


\subsection{Adaptive time stepping}
As can be seen from the previous examples, the phase-field models behavior can drastically change quickly over short timescales while varying very slowly over longer timescales. Various existing methods~\cite{CHENG201545,ETDRKEnergyStable} have utilized adaptive time stepping for faster simulation; we therefore implement daptive time stepping for the Cahn--Hilliard equation. We use the strategy motivated in in~\cite{Kostic2013}, which defines a dimensionless modified local truncation error and calculates the infinity norm of the local truncation error as a criterion. However, if their strategy were directly used in our simulation, division by zero would be (almost) inevitable and would lead to ineffective measurement. We therefore consider the further modified local truncation error
\begin{equation}
\label{eqn:lte}
\begin{split}
&\text{LTE}=\left|e_{k+1}-\frac{\tau}{\tau_{\text{old}}}e_k\right|,\\
\text{where}\quad e_{k+1}=\frac{\|\bm{u}^{k+1}-\bm{u}^k\|_\infty}{\|\bm{u}^k\|_\infty},\, &e_{k}=\frac{\|\bm{u}^{k}-\bm{u}^{k-1}\|_\infty}{\|\bm{u}^{k-1}\|_\infty}, \tau=t^{k+1}-t^k, \tau_{\text{old}}=t^{k}-t^{k-1}.
\end{split}
\end{equation}
Compared to the original formulation in~\cite{Kostic2013}, our formulation first applies the infinity norm to the finite difference errors $e_{k+1}$ and $e_{k}$, so that the division by zero issue is eliminated. 

Another noteworthy component of the adaptive time stepping strategy is the update of number of substeps $s$ for the SRK methods. Since the implicit method used in~\cite{Kostic2013} is unconditionally stable, their adaptive time stepping strategy does not need to account for the stability constraint. However, the stability constraints of our SRK methods are fully controlled by the number of substeps $s$, therefore we need to change $s$ each time the timestep $\tau$ increases or decreases to avoid redundant or insufficient number of substeps. Since the fourth-order term dominates in the Cahn--Hilliard, we choose our $s$ analytically according to the CFL condition given in~\cite{RKL2,SKARAS_RKG2021109879} with a small relaxation. The detailed method for changing $s$ is described in~\Cref{alg:substep_s}.

\begin{figure}[!t]
    \centering
    \begin{minipage}{0.8\linewidth}
        \begin{subfigure}{\linewidth}
        \captionsetup{justification=centering}
        \includegraphics[width=0.2\linewidth]{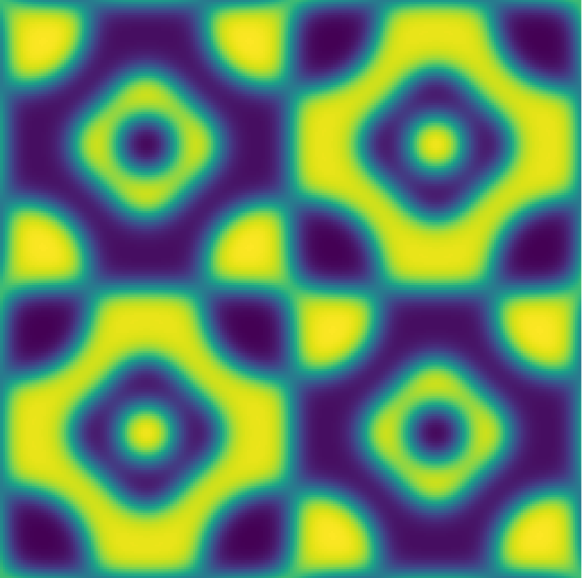}
        \hspace{0.7cm}
        \includegraphics[width=0.2\linewidth]{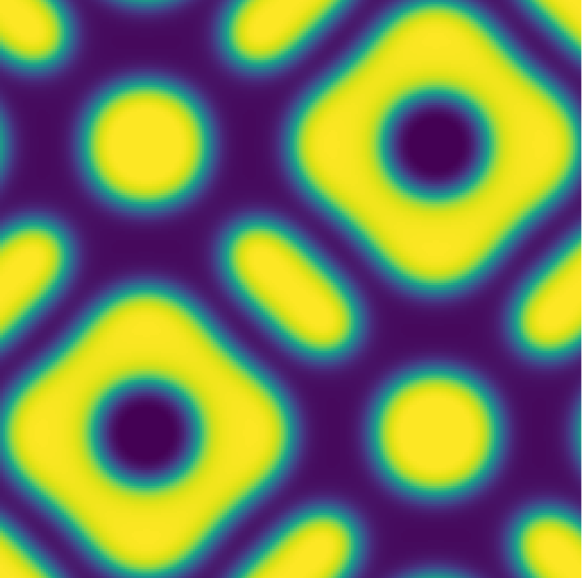}
        \hspace{0.7cm}
        \includegraphics[width=0.2\linewidth]{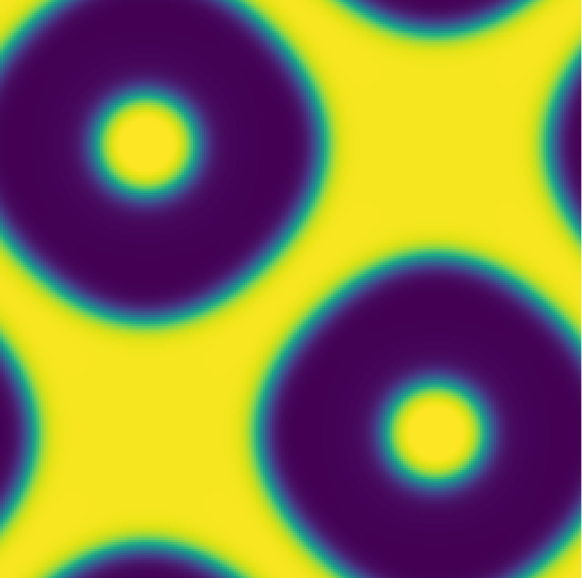}
        \hspace{0.7cm}
        \includegraphics[width=0.2\linewidth]{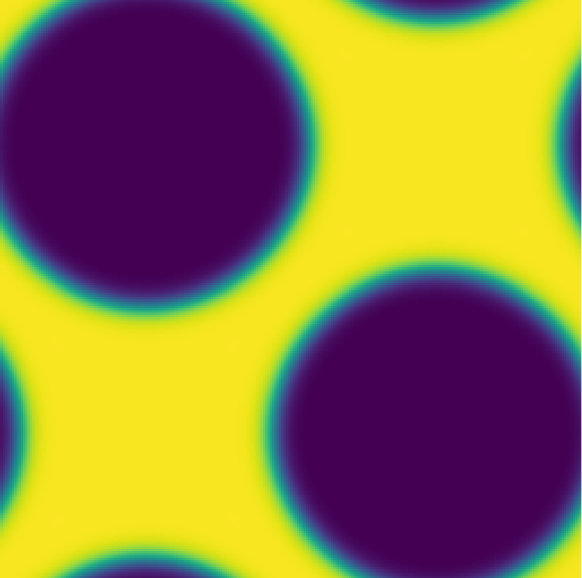}
        \caption{Cahn--Hilliard simulation adaptive time stepping, RKL2}
        \end{subfigure}
        \begin{subfigure}{\linewidth}
        \captionsetup{justification=centering}
        \includegraphics[width=0.2\linewidth]{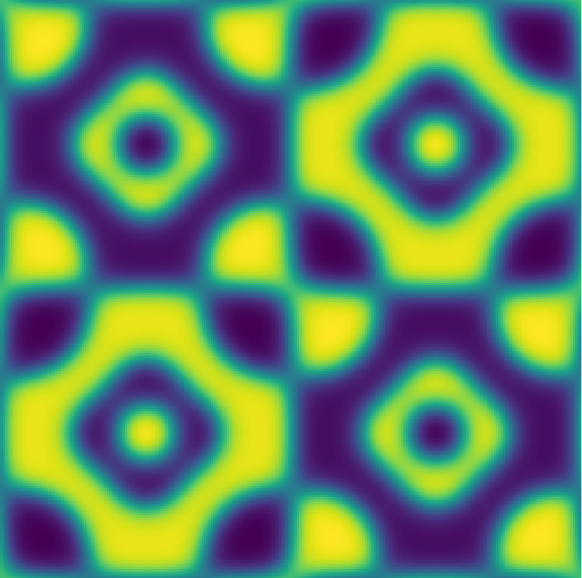}
        \hspace{0.7cm}
        \includegraphics[width=0.2\linewidth]{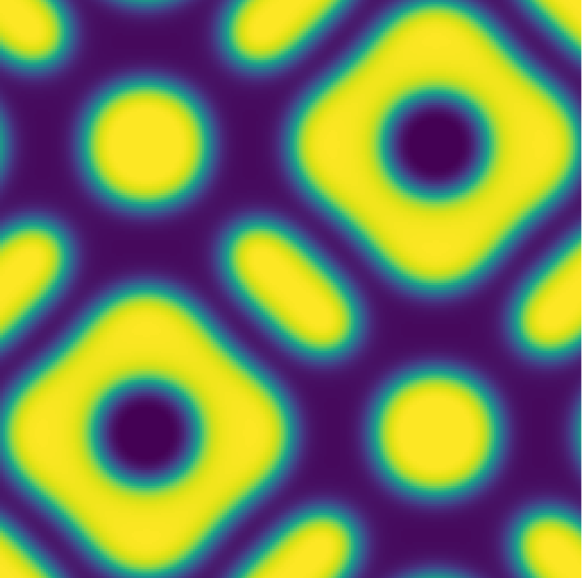}
        \hspace{0.7cm}
        \includegraphics[width=0.2\linewidth]{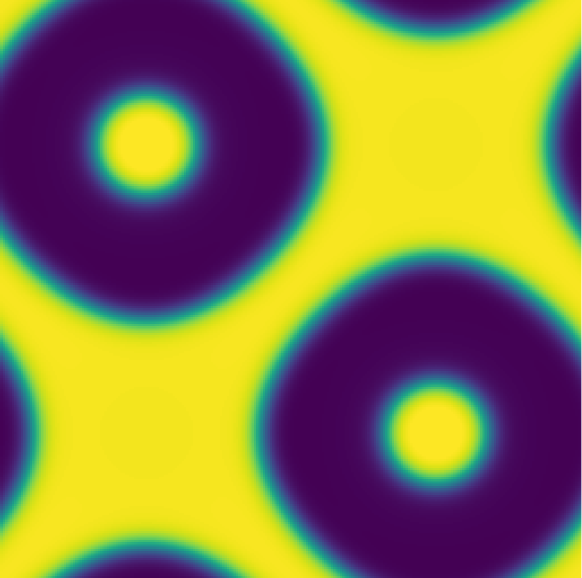}
        \hspace{0.7cm}
        \includegraphics[width=0.2\linewidth]{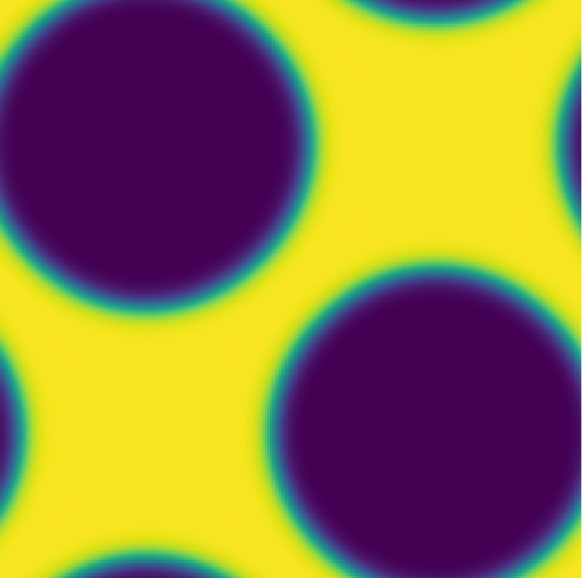}
        \caption{Cahn--Hilliard simulation adaptive time stepping, RKG2}
        \end{subfigure}
    \end{minipage}%
    \hspace{0.1cm}%
    \begin{minipage}{0.06\linewidth}
        \vspace{-0.5cm}
        \begin{tikzpicture}
            \begin{axis}[
                hide axis,
                scale only axis,
                height=6cm,
                width=0.1cm,
                colormap/viridis,
                colorbar right,
                colorbar style={
                    width=0.3cm,
                    ytick={-1, -0.5, 0, 0.5, 1},
                    ylabel={},
                },
                point meta min=-1,
                point meta max=1,
            ]
                \addplot [draw=none] coordinates {(0,0)};
            \end{axis}
        \end{tikzpicture}
    \end{minipage}
    \vspace{-1em}
    \caption{Cahn--Hilliard equation coarsening evolution at (approximately) $t=1, 2, 5$ and $8$, with initial condition~\Cref{eqn:sinusoidal_ic_ch} and adaptive time stepping. As shown, there is no visual difference between the RKL2 and RKG2 simulations, and the results are similar to the constant timestep simulations in~\Cref{fig:ch_sinusoidal}.}
    \label{fig:ch_sinusoidal_adaptive}
\end{figure}
\begin{figure}[!t]
    \centering
    \begin{subfigure}{0.4\linewidth}
    \captionsetup{justification=centering}
    \begin{overpic}[width=1\textwidth]{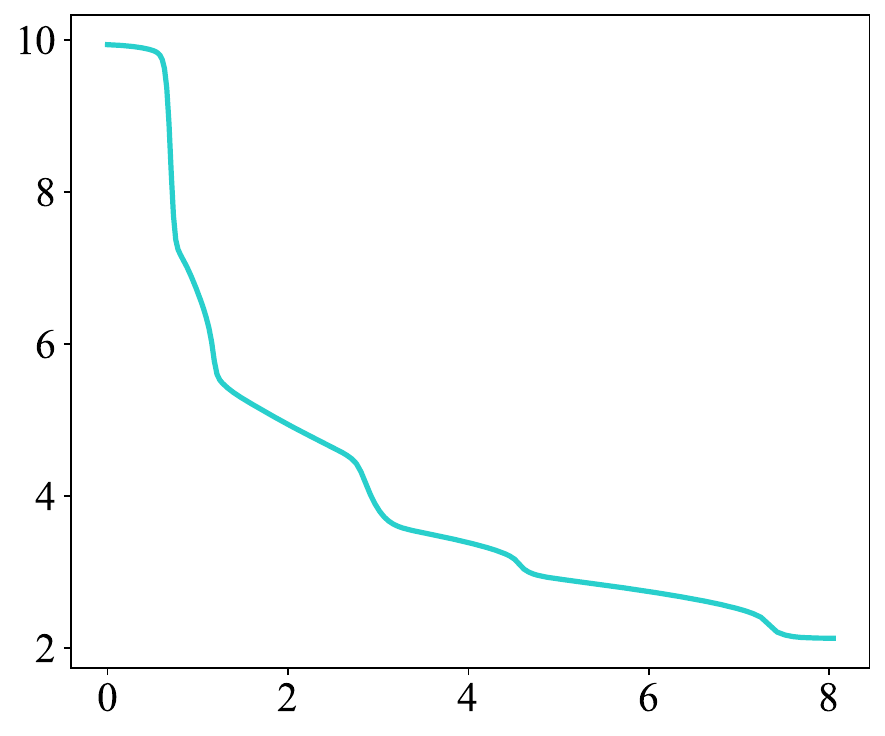}
        \put(52,-2){\small $t$}
        \put(-4,38){\rotatebox{90}{\small Energy}}
    \end{overpic}
    \caption{Cahn--Hilliard equation energy, RKL2}
    \end{subfigure}
    \hspace{0.5em}
    \begin{subfigure}{0.4\linewidth}
    \captionsetup{justification=centering}
    \begin{overpic}[width=1\textwidth]{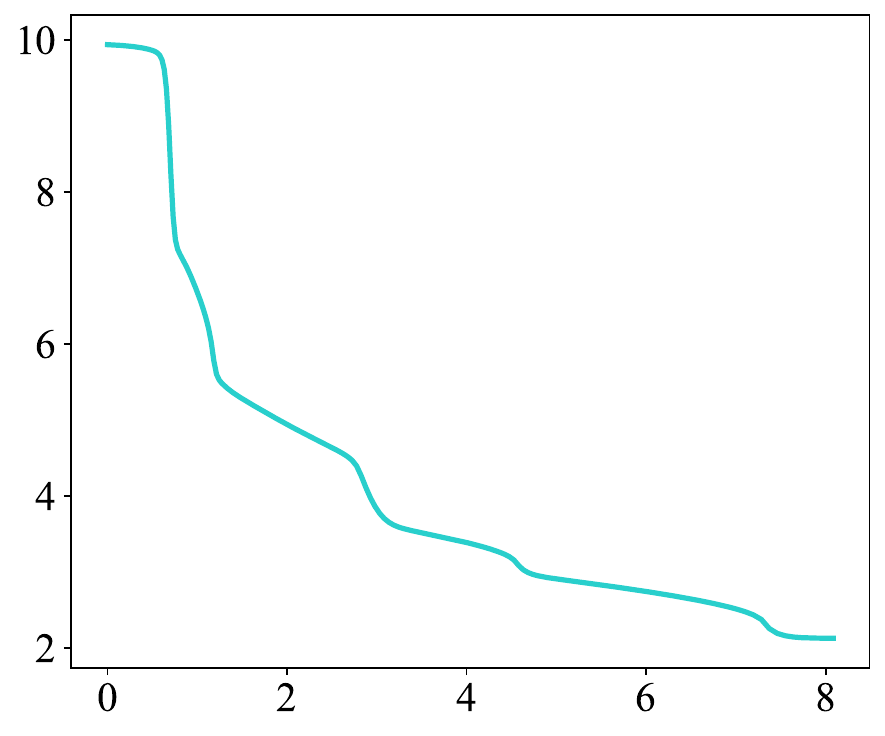}
        \put(52,-2){\small $t$}
        \put(-4,38){\rotatebox{90}{\small Energy}}
    \end{overpic}
    \caption{Cahn--Hilliard equation energy, RKG2}
    \end{subfigure}
    \caption{Energy dissipation of the Cahn--Hilliard equation solutions with adaptive time stepping.}
    \label{fig:ch_sinusoidal_adaptive_energy}
\end{figure}
\begin{algorithm}[!t]
\caption{Substep $s$ changing procedure}\label{alg:substep_s}
\begin{algorithmic}
\STATE \textbf{Require:} Timestep $\tau$ and spatial grid size $h$; old substep number $s$, method $RKL$ or $RKG$; case $\max$ or $\min$.
\STATE \textbf{Step 1:} Compute CFL number $C=\epsilon^2\frac{\tau}{h^4}$.
\STATE \textbf{Step 2:} If method is $RKL$, set $s_{\text{new}}=\text{ceil}\left(\frac{-1+\sqrt{9+512C}}{2}\right)$. If method is $RKG$, set $s_{\text{new}}=\text{ceil}\left(\frac{-3+\sqrt{25+768C}}{2}\right)$.
\STATE \textbf{Step 3:} Relax the condition a little bit by setting $s_{\text{new}}=\text{ceil}(1.02s)$.
\STATE If case is $\max$, set $s=\max\{s_{\text{new}}, s\}$. If case is $\min$, set $s=\min\{s_{\text{new}}, s\}$.
\end{algorithmic}
\end{algorithm}

\begin{figure}[!t]
    \centering
    \begin{subfigure}{0.4\linewidth}
    \captionsetup{justification=centering}
    \begin{overpic}[width=1\textwidth]{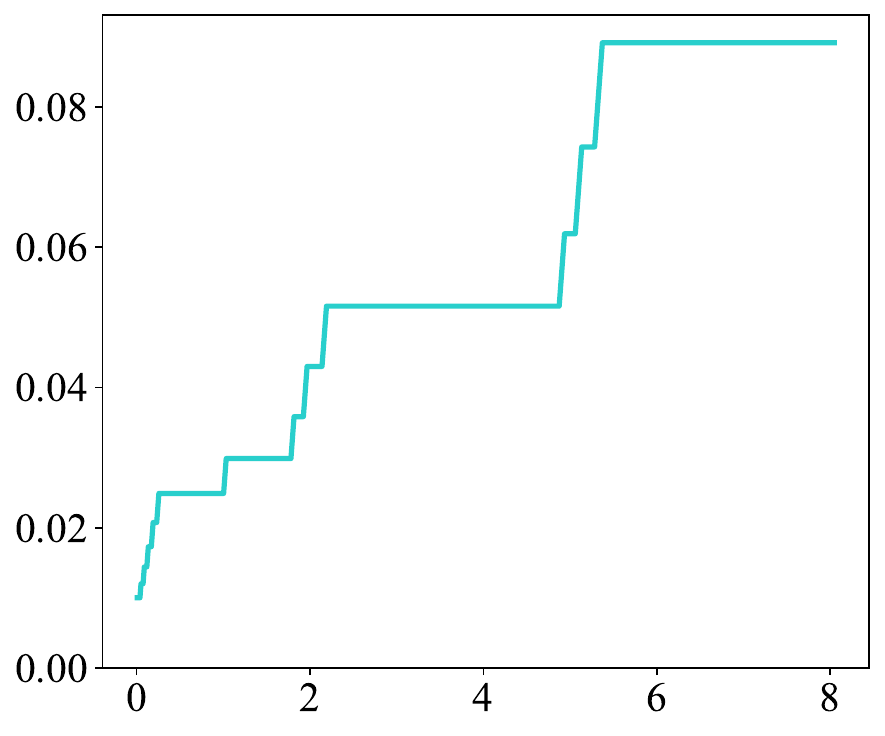}
        \put(52,-2){\small $t$}
        \put(-4,38){\rotatebox{90}{\small Step size}}
    \end{overpic}
    \caption{Change of timestep size for the Cahn--Hilliard equation, RKL2}
    \end{subfigure}
    \hspace{0.5em}
    \begin{subfigure}{0.4\linewidth}
    \captionsetup{justification=centering}
    \begin{overpic}[width=1\textwidth]{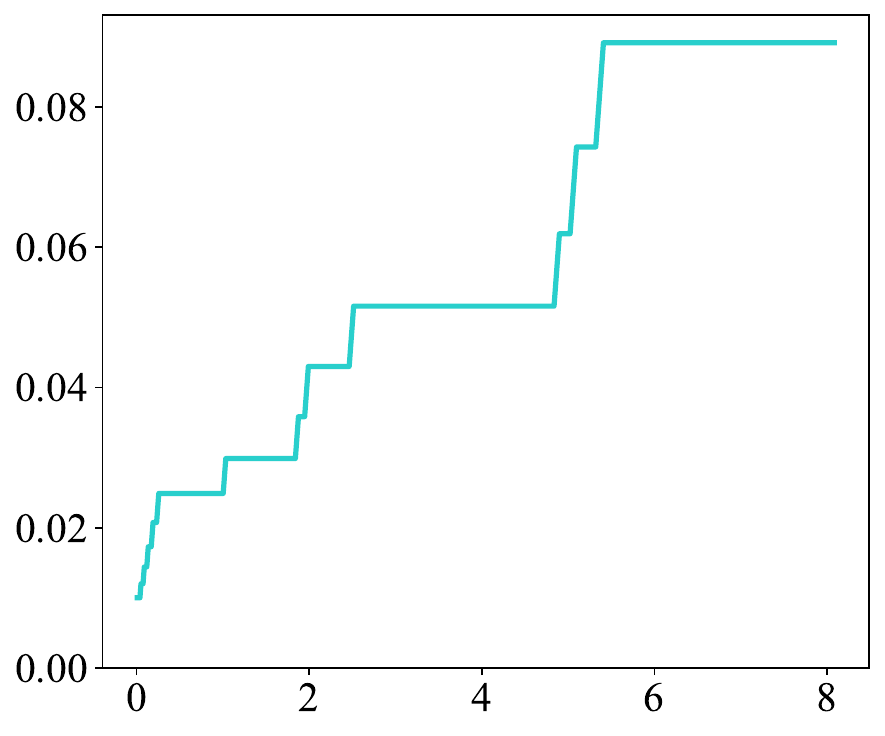}
        \put(52,-2){\small $t$}
        \put(-4,38){\rotatebox{90}{\small Step size}}
    \end{overpic}
    \caption{Change of timestep size for the Cahn--Hilliard equation, RKG2}
    \end{subfigure}
    \caption{Evolution of $\tau$ over the adaptive time stepping simulation of the Cahn--Hilliard equation, with sinusoidal initial condition~\Cref{eqn:sinusoidal_ic_ch}.}
    \label{fig:ch_sinusoidal_adaptive_timestep}
\end{figure}

We summarize the algorithm described above in~\Cref{alg:adaptive_time_stepping}, and conduct Cahn--Hilliard equation simulations accordingly. We choose the adaptive timestep tolerance \textit{tol} to be $1e$-$3$ and \textit{MaxCount} to be 3. We again choose the initial condition to be the perturbed sinusoidal initial condition~\Cref{eqn:sinusoidal_ic_ch}, with the same initial parameters as in the constant timestep case. As shown in~\Cref{fig:ch_sinusoidal_adaptive} and~\Cref{fig:ch_sinusoidal_adaptive_energy}, the adaptive time stepping simulations and the energy curves match the constant timestep results in~\Cref{fig:ch_sinusoidal} and~\Cref{fig:ch_sinusoidal_energy}. Moreover, we plot the $\tau$ used for each timestep for both methods in~\Cref{fig:ch_sinusoidal_adaptive_timestep}; as shown in the figures, the $\tau$ used significantly increases from $T=0$ to $T=8$. We also conduct a wall-clock comparison of the simulation times. For RKL2 the constant timestep version uses 896.02s, while the adaptive timestep version uses 344.95s. For RKG2, the constant timestep version uses 1185.29s, while the adaptive timestep version uses 502.32s. This shows the effectiveness of our adaptive timestepping strategy in the Cahn--Hilliard equation simulations.

\begin{algorithm}[!t]
\caption{Adaptive time stepping procedure}\label{alg:adaptive_time_stepping}
\begin{algorithmic}
\STATE \textbf{Require:} Numerical solutions of two steps $\bm{u}^{k},\bm{u}^{k-1}$; current and previous timesteps $\tau, \tau_{\text{old}}$; adaptive time tolerance $tol$; small error counting $count$; maximum number of count $MaxCount$.
\STATE \textbf{Step 1:} Compute $\bm{u}^{k+1}$ using timestep $\tau$.
\STATE \textbf{Step 2:} Calculate LTE according to~\Cref{eqn:lte}.
\STATE \textbf{Step 3:} If the original energy is not monotone, set $\tau_{\text{old}}=\tau$ and $\tau=0.8\tau$. Set $s$ with $\tau$ and case min, then go back to step 1.
\STATE \textbf{Step 4:} Else, if $\text{LTE}<tol$, $count=count+1$. If $count$ reaches $MaxCount$, set $\tau_{\text{old}}=\tau$, $\tau=1.2\tau$, and $count=0$.
\STATE \textbf{Step 5:} Set $s$ with $\tau$ and case max, then go back to step 1.
\STATE \textbf{Step 6:} Endif.
\end{algorithmic}
\end{algorithm}

\section{Conclusion}
\label{sec:conclusion}
In this paper, we present simulations of Allen--Cahn and Cahn--Hilliard equations using the explicit stabilized Runge--Kutta methods. Theoretical analysis guaranties that the discrete maximum principle is preserved for the Allen--Cahn equation with Strang splitting within the stability constraint. We also prove that the modified discrete energy, under two different ways of splitting, dissipate with a minor tradeoff of maximum timestep allowed, by estimating the positive-definiteness of the stability polynomial. Numerical simulations also verify and validate the rate of convergence and the energy dissipation of both Allen--Cahn and Cahn--Hilliard equations, for both rectangular domains and irregular domains, under periodic and Neumann boundary conditions.

For future work, one can potentially develop higher-order stabilized Runge--Kutta methods specifically designed for gradient flows. Another direction can be to prove dissipation of the original energy for Allen--Cahn and the Cahn--Hilliard equations, or find out the limitation of the methods, since only empirical evidence is provided for 2D maximum principle and Cahn--Hilliard energy dissipation. Moreover, our methods can be applied to higher-order parabolic equations with variable coefficients easily, and gain more efficiency when compared with the implicit methods. Finally, theoretical guarantee of the maximum principle when applying RKG methods to the Allen--Cahn equation with Dirichlet boundary condition can also help us better understand the stabilized Runge--Kutta methods, and it would also be beneficial if the maximum principle for RKL and RKG methods can be proved for more general domain shapes.


\section*{Acknowledgment} 
\noindent
The work of the second author was supported by the Advanced Simulation and Computing Program (ASC) at Los Alamos National Laboratory, operated by Triad National Security, LLC, for the National Nuclear Security Administration of the U.S. Department of Energy contract 89233218CNA000001. The work of the first and the third author was supported by ONR grant DOD N00014-23-1-2565.

The first author used AI-assisted tools, including ChatGPT and Claude, for support with function visualizations in theorem proving, English language polishing, LaTeX editing, and partial code checking. 
The third author used Claude for one stage of proofreading but not the final proofreading of the paper, which was checked by a human editor.

\newpage

\bibliographystyle{elsarticle-num} 
\bibliography{reference}

\end{document}